\documentclass[11pt]{article}
\usepackage{graphicx}
\usepackage{amsmath}
\usepackage{amsfonts}
\usepackage{enumerate}
\usepackage{latexsym}
\usepackage{epsfig}

\newtheorem{theorem}{Theorem}[section]
\newtheorem{lemma}[theorem]{Lemma}
\newtheorem{coroll}[theorem]{Corollary}

\def\proofbox{\begin{picture}(6.5,6.5)
\put(0,0){\framebox(6.5,6.5){}}\end{picture}}
\newenvironment{proof}{\noindent{\it Proof.\quad}}{\hfill\proofbox}

\begin{document}
\begin{center}

{\bf \Large Superinjective Simplicial Maps of Complexes of}

\vspace{0.1in}
{\bf \Large Curves and Injective Homomorphisms of}\\
\vspace{0.1in}
{\bf \Large Subgroups of Mapping Class Groups II }\\
\vspace{0.2in}

\large {Elmas Irmak} \\
\end{center}
\vspace{0.1in}



\begin{abstract}
Let $R$ be a compact, connected, orientable surface of genus $g$
with $p$ boundary components. Let $\mathcal{C}(R)$ be the complex
of curves on $R$ and $Mod_R^*$ be the extended mapping class group
of $R$. Suppose that either $g = 2$ and $p \geq 2$ or $g \geq 3$
and $p \geq 0$. We prove that a simplicial map $\lambda:
\mathcal{C}(R) \rightarrow \mathcal{C}(R)$ is superinjective if
and only if it is induced by a homeomorphism of $R$. As a
corollary, we prove that if $K$ is a finite index subgroup of
$Mod_R^*$ and $f : K \rightarrow Mod_R^*$ is an injective
homomorphism, then $f$ is induced by a homeomorphism of $R$ and
$f$ has a unique extension to an automorphism of $Mod_R^*$. This
extends the author's previous results about closed connected
orientable surfaces of genus at least 3, to the surface $R$.
\end{abstract}
\vspace{0.12in}

{\it Mathematics Subject Classification}: 57M99, 20F38.\\

\section{Introduction}

Let $R$ be a compact, connected, orientable surface of genus $g$
with $p$ boundary components. The mapping class group, $Mod_R$, of
$R$ is the group of isotopy classes of orientation preserving
homeomorphisms of $R$. The pure mapping class group, $PMod_R$, is
the subgroup of $Mod_R$ consisting of isotopy classes of
homeomorphisms which preserve each boundary component of $R$. The
extended mapping class group, $Mod_R^*$, of $R$ is the group of
isotopy classes of all
(including orientation reversing) homeomorphisms of $R$.\\

Assume that either $g = 2$ and $p \geq 2$ or $g \geq 3$ and $p
\geq 0$. The main results of the paper are the following two
theorems:

\begin{theorem} \label{theorem1} A simplicial map, $\lambda : \mathcal{C}(R)
\rightarrow \mathcal{C}(R)$, is superinjective if and only if
$\lambda$ is induced by a homeomorphism of $R$.
\end{theorem}

\begin{theorem} \label{theorem2} Let $K$ be a finite index subgroup of
$Mod_R^*$ and $f$ be an injective homomorphism, $f:K \rightarrow
Mod_R^*$. Then $f$ is induced by a homeomorphism of the surface
$R$ (i.e. for some $g \in Mod_R^*$, $f(k) = gkg^{-1}$ for all $k
\in K$) and $f$ has a unique extension to an automorphism of
$Mod_R^*$.\end{theorem}

These two theorems were proven for closed, connected, orientable
surfaces of genus at least 3 by the author in \cite{Ir}. They were
motivated by the following theorems of Ivanov and the theorem of
Ivanov and McCarthy.

\begin{theorem} {\bf (Ivanov) \cite{Iv1}} \label{theoremA} Let $R$
be a compact, orientable surface possibly
with nonempty boundary. Suppose that the genus of $R$ is at least
2. Then, all automorphisms of $ \mathcal{C}(R)$ are given by
elements of $Mod_R^*$. More precisely, if $R$ is not a closed
surface of genus 2, then $Aut(\mathcal{C}(R))= Mod_R ^*$. If $R$
is a closed surface of genus 2, then $Aut(\mathcal{C}(R))= Mod_R
^* /Center(Mod_R ^*)$.
\end{theorem}

\begin{theorem} {\bf (Ivanov) \cite{Iv1}} \label{theoremB} Let $R$
be a compact, orientable surface possibly with nonempty boundary.
Suppose that the genus of $R$ is at least 2 and $R$ is not a
closed surface of genus 2. Let $\Gamma_1, \Gamma_2$ be finite
index subgroups of $Mod_R ^*$. Then, all isomorphisms $\Gamma_1
\rightarrow \Gamma_2$ have the form $x \rightarrow gxg^{-1}, g \in
Mod_R^*$.
\end{theorem}

\begin{theorem} {\bf (Ivanov, McCarthy) \cite{IMc}} \label{theoremC}
Let $R$ and $R'$ be compact, connected, orientable surfaces.
Suppose that the genus of $R$ is at least 2, $R'$ is not a closed
surface of genus 2, and the maxima of ranks of abelian subgroups
of $Mod_R$ and $Mod_{R'}$ differ by at most one. If $h: Mod_R
\rightarrow Mod_{R'}$ is an injective homomorphism, then  $h$ is
induced by a homeomorphism $H: S \rightarrow R'$, (i.e.
$h([G])=[HGH^{-1}]$, for every orientation preserving
homeomorphism $G : R \rightarrow R$).
\end{theorem}

For the surfaces that we consider in this paper, Theorem
\ref{theorem1} generalizes Ivanov's Theorem \ref{theoremA},
Theorem \ref{theorem2} generalizes Ivanov's Theorem \ref{theoremB}
and Ivanov and McCarthy's Theorem \ref{theoremC} (in the case when
the surfaces are the same).\\

\indent In section 2, we give some properties of superinjective
simplicial maps of the complex of curves, $\mathcal{C}(R)$.\\
\indent In section 3, we prove that a superinjective simplicial
map $\lambda : \mathcal{C}(R) \rightarrow \mathcal{C}(R)$, induces
an injective simplicial map on the complex of arcs,
$\mathcal{B}(R)$, and by using this map we prove that $\lambda$ is
induced by a homeomorphism of $R$.\\
\indent In section 4, we prove that if $K$ is a finite index
subgroup of $Mod_R^*$ and $f:K \rightarrow Mod_R^*$ is an
injective homomorphism, then $f$ induces a superinjective
simplicial map of $\mathcal{C}(R)$, and by using this map we prove
that $f$ is induced by a homeomorphism of $R$.

\section{Superinjective Simplicial Maps of Complexes of Curves}

Let $R$ be a compact, connected, orientable surface of genus $g$
with $p$ boundary components. We assume that either $g = 2$ and $p
\geq 2$ or $g \geq 3$ and $p \geq 0$. Since our main results have
already been proven for $g \geq 3$ and $p = 0$ in \cite{Ir}, we
will assume that $p > 0$ when $g \geq 3$ throughout the paper.\\

\indent A \textit{circle} on $R$ is a properly embedded image of
an embedding $S^{1} \rightarrow R$. A circle on $R$ is said to be
\textit{nontrivial} (or \textit{essential}) if it does not bound a
disk and it is not homotopic to a boundary component of $R$. Let
$C$ be a collection of pairwise disjoint circles on $R$. The
surface obtained from $R$ by cutting along $C$ is denoted by
$R_C$. A nontrivial circle $a$ on $R$ is called
\textit{(k,m)-separating} (or a \textit{(k,m) circle}) if the
surface $R_{a}$ is disconnected and one of its components is a
genus $k$ surface with $m$ boundary components, where $ 1 \leq k
\leq g$, $ 1 \leq m < p$. If $R_a$ is connected, then $a$ is
called \textit{nonseparating}. Let $\mathcal{A}$ denote the set of
isotopy classes of nontrivial circles on $R$. The
\textit{geometric intersection number} $i(\alpha, \beta)$ of
$\alpha, \beta \in \mathcal{A}$ is the minimum number of points of
$a \cap b$ where $a \in \alpha$ and $b \in \beta$.\\
\indent The \textit{complex of curves}, $\mathcal{C}(R)$, on $R$,
introduced by W. Harvey \cite{Ha}, is an abstract simplicial
complex with vertex set $\mathcal{A}$ such that a set of $n$
vertices $\{{ \alpha_{1}}, {\alpha_{2}}, ..., {\alpha_{n}}\}$
forms an $n-1$ simplex if and only if ${\alpha_{1}},
{\alpha_{2}},..., {\alpha_{n}}$ have pairwise disjoint
representatives.\\

\noindent {\bf Definition:} A simplicial map $\lambda :
\mathcal{C}(R) \rightarrow \mathcal{C}(R)$ is called {\bf
superinjective} if the following condition holds: if $\alpha,
\beta$ are two vertices in $\mathcal{C}(R)$ such that
$i(\alpha,\beta) \neq 0$, then $i(\lambda(\alpha),\lambda(\beta))
\neq 0$.\\

In this section, we show some properties of superinjective
simplicial maps of $ \mathcal{C}(R)$. The proofs of the Lemmas
2.1-2.4 are similar to the proofs of the corresponding lemmas in
the closed case which are given in \cite{Ir}. So, we will only
state them here.

\begin{lemma}
\label{injective} A superinjective simplicial map, $\lambda :
\mathcal{C}(R) \rightarrow \mathcal{C}(R)$, is injective.
\end{lemma}

\begin{lemma}
\label{conbyanedge} Let $\alpha, \beta$ be two distinct vertices
of $\mathcal{C}(R)$, and let $\lambda : \mathcal{C}(R) \rightarrow
\mathcal{C}(R)$ be a superinjective simplicial map. Then, $\alpha$
and $\beta$ are connected by an edge in $\mathcal{C}(R)$ if and
only if $\lambda(\alpha)$ and $\lambda(\beta)$ are connected by an
edge in $\mathcal{C}(R)$.
\end{lemma}

Let $P$ be a set of pairwise disjoint circles on $R$. $P$ is
called a {\it pair of pants decomposition} of $R$, if $R_P$ is a
disjoint union of genus zero surfaces with three boundary
components, pairs of pants. A pair of pants of a pants
decomposition is the image of one of these connected components
under the quotient map $q:R_P \rightarrow R$ together with the
image of the boundary components of this component. The image of
the boundary of this component is called the \textit{boundary of
the pair of pants}. A pair of pants is called \textit{embedded} if
the restriction of $q$ to the corresponding component of $R_P$ is
an embedding. An ordered set $(a_1, ..., a_{3g-3+p})$ is called an
{\it ordered pair of pants decomposition} of $R$ if $\{a_1, ...,
a_{3g-3+p}\}$ is a pair of pants decomposition of $R$.

\begin{lemma}
\label{imageofpantsdecomp} Let $\lambda : \mathcal{C}(R)
\rightarrow \mathcal{C}(R)$ be a superinjective simplicial map.
Let $P$ be a pair of pants decomposition of $R$. Then, $\lambda$
maps the set of isotopy classes of elements of $P$ to the set of
isotopy classes of elements of a pair of pants decomposition,
$P'$, of $R$.
\end{lemma}

Let $P$ be a pair of pants decomposition of $R$. Let $a$ and $b$
be two distinct elements in $P$. Then, $a$ is called {\it
adjacent} to $b$ w.r.t. $P$ iff there exists a pair of pants in
$P$ which has $a$ and $b$ on its boundary.\\

\noindent {\bf Remark:} Let $P$ be a pair of pants decomposition
of $R$. Let $[P]$ be the set of isotopy classes of elements of
$P$. Let $\alpha, \beta \in [P]$. We say that $\alpha$ is adjacent
to $\beta$ w.r.t. $[P]$ if the representatives of $\alpha$ and
$\beta$ in $P$ are adjacent w.r.t. $P$. By Lemma
\ref{imageofpantsdecomp}, $\lambda$ gives a correspondence on the
isotopy classes of elements of pair of pants decompositions of
$R$. $\lambda([P])$ is the set of isotopy classes of elements of a
pair of pants decomposition which corresponds to $P$, under this
correspondence.

\begin{lemma}
\label{adjacent} Let $\lambda : \mathcal{C}(R) \rightarrow
\mathcal{C}(R)$ be a superinjective simplicial map. Let $P$ be a
pair of pants decomposition of $R$. Then, $\lambda$ preserves the
adjacency relation for two circles in $P$, i.e. if $a, b \in P$
are two circles which are adjacent w.r.t. $P$ and $[a]=\alpha,
[b]=\beta$, then $\lambda(\alpha), \lambda(\beta)$ are adjacent
w.r.t. $\lambda([P])$.\end{lemma}

Let $P$ be a pair of pants decomposition of $R$. A curve $x \in P$
is called a {\it 4-curve} in $P$, if there exist four distinct
circles in $P$, which are adjacent to $x$ w.r.t. $P$.

\begin{lemma}
\label{division} Let $\lambda : \mathcal{C}(R) \rightarrow
\mathcal{C}(R)$ be a superinjective simplicial map. Let $a$ be a
$(k,1)$-separating circle on $R$, where $2 \leq k \leq g$. Let
$R_1, R_2$ be the subsurfaces of $R$ s.t. $R_1$ has genus $k$ and
has $a$ as its boundary, and $R_2 = R \setminus R_1 \cup a$. Let
$a' \in \lambda([a])$. Then $a'$ is a $(k,1)$-separating circle
and there exist subsurfaces $R_1', R_2'$ of $R$ s.t. $R_1'$ has
genus $k$ and has $a'$ as its boundary, $R_2' = R \setminus R_1'
\cup a'$, $\lambda(\mathcal{C}(R_1)) \subseteq \mathcal{C}(R_1')$
and $\lambda(\mathcal{C}(R_2)) \subseteq \mathcal{C}(R_2')$.
\end{lemma}

\begin{proof} Let $a$ be a $(k,1)$-separating circle on $R$ and
$2 \leq k \leq g$. Let $R_1$ be the subsurface of $R$ of genus $k$
which has $a$ as its boundary and $R_2 = R \setminus R_1 \cup a$.
If $k=2$, we choose a pair of pants decomposition $P_1 = \{a_1,
a_2, a_3, a_4 \}$ of $R_1$ as shown in Figure 1, (i) and then
complete $P_1 \cup \{a\}$ to a pair of pants decomposition $P$ of
$R$ in any way we like. If $k \geq 3$, we choose a pair of pants
decomposition $P_1 = \{a_1, a_2, ..., a_{3k-2} \}$ of $R_1$ and
then complete $P_1 \cup \{a\}$ to a pair of pants decomposition
$P$ of $R$ such that each of $a_i$ is a 4-curve in $P$ for
$i=1,2,...,3k-2$ and $a, a_1, a_3$ are the boundary components of
a pair of pants of $P_1$. In Figure 1 (ii), we show how to choose
$P_1$ when $k=4$. In the other cases, when $k=3$ or $k \geq 5$, a
similar pair of pants decomposition of $R_1$ can be chosen.\\

\begin{figure}[htb]
\label{a1 }
\begin{center}
\epsfxsize=1.85in \epsfbox{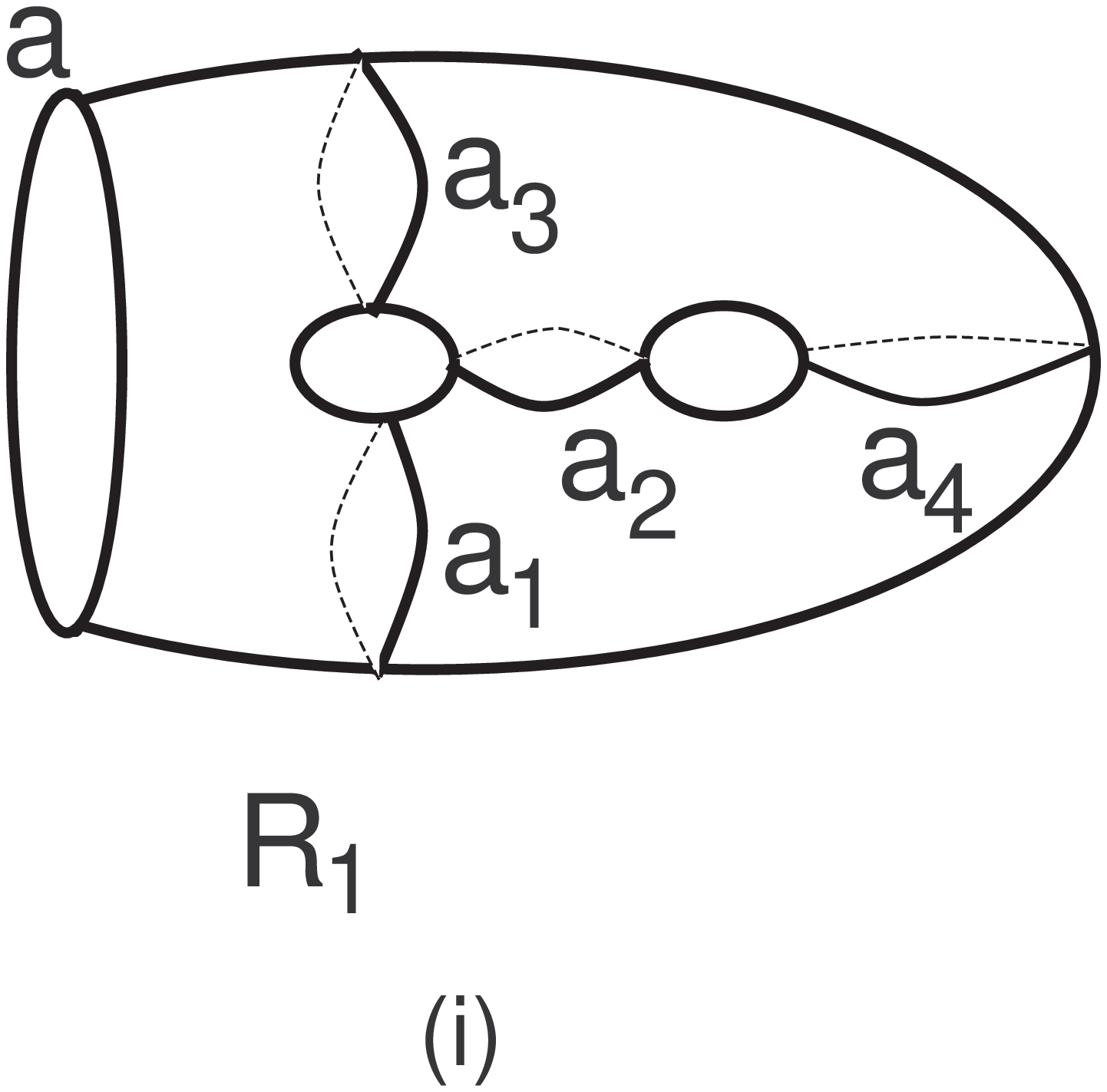} \epsfxsize=2.85in
\epsfbox{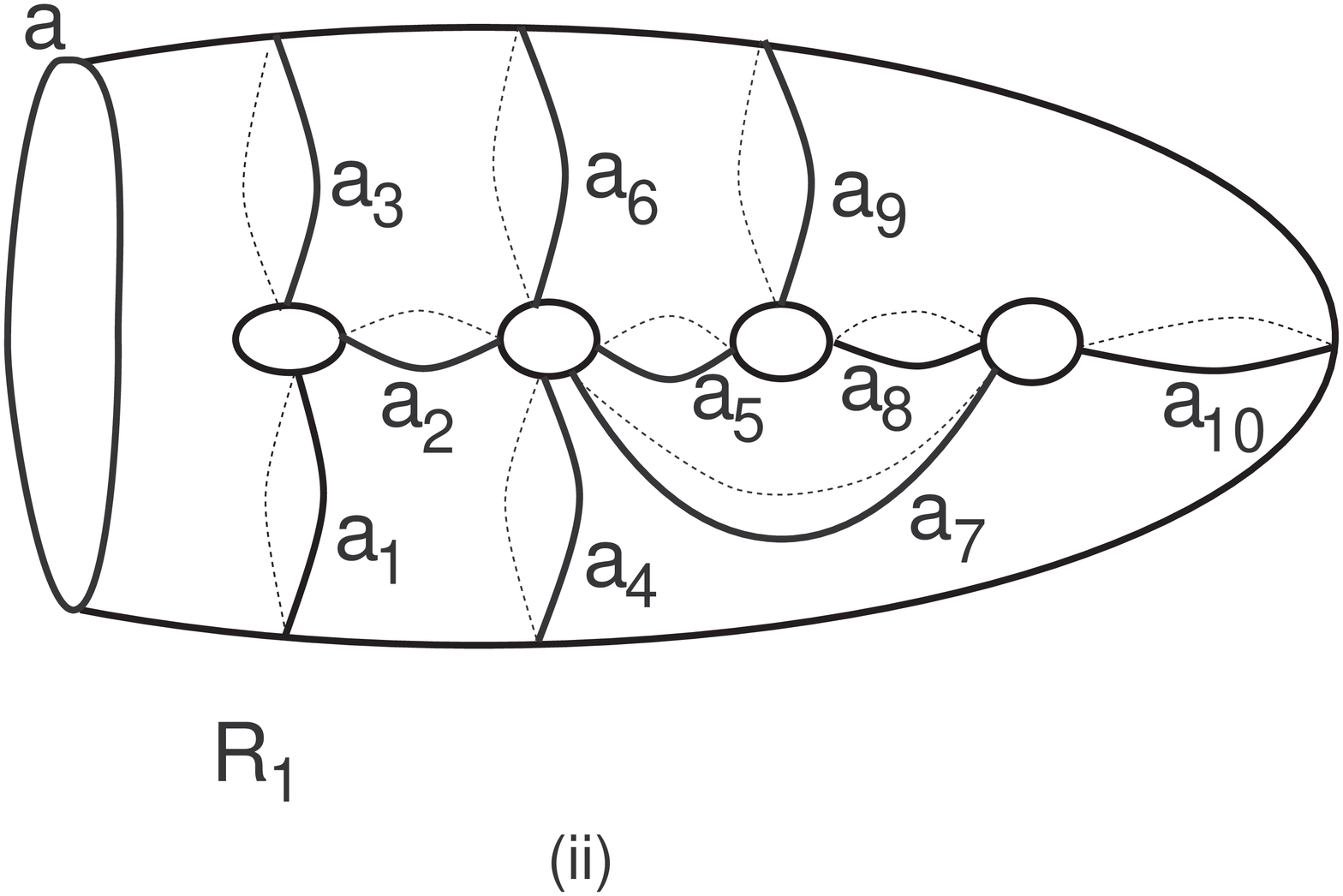}

\caption{Pants decompositions}
\end{center}
\end{figure}

Let $P'$ be a pair of pants decomposition of $R$ such that
$\lambda([P]) = [P']$. Let $a_i'$ be the representative of
$\lambda([a_i])$ which is in $P'$, for $i=1,..., 3k-2$, and $a'$
be the representative of $\lambda([a])$ which is in $P'$. Let
$P_1'= \{a_1', a_2', ..., a_{3k-2}'\}.$ By using Lemma 2.4 and
following the proof of Lemma 3.5 in [3], we can see that $a'$ is a
$(k,1)$-separating circle, there exist subsurfaces $R_1', R_2'$ of
$R$ s.t. $R_1'$ has genus $k$ and has $a'$ as its boundary, $R_2'
= R \setminus R_1' \cup a'$ and $P_1'$ is a pair of pants
decomposition of $R_1'$. Let $P_2 = P \setminus (P_1 \cup {a})$
and $P_2' = P' \setminus (P_1' \cup {a'})$. Then $P_2, P_2'$ are
pair of pants decompositions of $R_2, R_2'$ respectively as $P_1,
P_1'$ are pair of pants
decompositions of $R_1, R_1'$ respectively.\\

Now, let $\alpha$ be a vertex in $\mathcal{C}(R_1)$. Then, either
$\alpha \in [P_1]$ or $\alpha$ has a nonzero geometric
intersection with an element of $[P_1]$. In the first case,
clearly $\lambda(\alpha) \in \mathcal{C}(R_1')$ since elements of
$[P_1]$ correspond to elements of $[P_1'] \subseteq
\mathcal{C}(R_1')$. In the second case, since $\lambda$ preserves
zero and nonzero geometric intersection (since $\lambda$ is
superinjective) and $\alpha$ has zero geometric intersection with
the elements of $[P_2]$ and $[a]$, and nonzero intersection with
an element of $[P_1]$, $\lambda(\alpha)$ has zero geometric
intersection with elements of $[P_2']$ and $[a']$, and nonzero
intersection with an element of $[P_1']$. Then, $\lambda(\alpha)
\in \mathcal{C}(R_1')$. Hence, $\lambda(\mathcal{C}(R_1))
\subseteq \mathcal{C}(R_1')$. The proof of
$\lambda(\mathcal{C}(R_2)) \subseteq \mathcal{C}(R_2')$ is
similar.\end{proof}

\begin{lemma}
\label{separating1} Let $\lambda : \mathcal{C}(R) \rightarrow
\mathcal{C}(R)$ be a superinjective simplicial map. Then,
$\lambda$ sends the isotopy class of a nonseparating circle to the
isotopy class of a nonseparating circle.\end{lemma}

\begin{proof} Let $c$ be a nonseparating circle. Let's choose a
separating (2, 1) circle, $a$ on $R$ s.t. $c$ is in genus 2
subsurface, $R_1$, bounded by $a$. Then, $c$ can be completed to a
pants decomposition $P_1=\{a_1, a_2, a_3, a_4\}$ on $R_1$, where
$a_4 = c$ as in Figure 1, (i). Then we can complete $P_1 \cup
\{a\}$ to a pair of pants decomposition, $P$, on $R$ in any way we
like. Let $P'$ be a pair of pants decomposition of $R$ such that
$\lambda([P]) = [P']$. Let $a_i'$ be the representative of
$\lambda([a_i])$ which is in $P'$, for $i=1,..., 4$, and $a'$ be
the representative of $\lambda([a])$ which is in $P'$. Let $P_1'=
\{a_1', a_2', a_3', a_4'\}.$ By Lemma 2.5, $a'$ is a (2,1)
separating circle bounding a subsurface, say $R_1'$, and $P_1'$ is
a pants decomposition on $R_1'$. Then, by using the adjacency
relation between the elements of $P_1' \cup \{a'\}$ it is easy to
see that $a_4' \in \lambda([c])$ is a nonseparating circle.\end{proof}\\

For every 4-curve $x$ in a pants decomposition $P$, there exist
two pairs of pants $A(x)$ and $B(x)$ of $P$ having $x$ as one of
their boundary components. Let $C(x)= A(x) \cup B(x)$. The
boundary of $C(x)$ consists of four distinct curves, which are
adjacent to $x$ w.r.t. $P$.

\begin{lemma}
\label{ ?} Let $\lambda : \mathcal{C}(R) \rightarrow
\mathcal{C}(R)$ be a superinjective simplicial map. Then,
$\lambda$ sends the isotopy class of a $(0,3)$-separating circle
to the isotopy class of a $(0,3)$-separating circle.\end{lemma}

\begin{proof} Let $c$ be a $(0,3)$-separating circle. A $(0,3)$-separating
circle exists if the number of boundary components, $p$, of $R$,
is at least 2. So, we assume that $p \geq 2$. If $p=2$, then it is
easy to see that a circle is a $(0,3)$-separating circle iff it is
a $(g,1)$ separating circle on $R$. Since $c$ is a $(0,3)$
separating circle, $c$ is a $(g,1)$ separating circle on $R$.
Then, by Lemma \ref{division}, $\lambda([c])$ has a $(g,1)$
representative which is a $(0,3)$-separating circle on $R$.\\

\begin{figure}[htb]
\label{a0 }
\begin{center}
\epsfxsize=3.5in \epsfbox{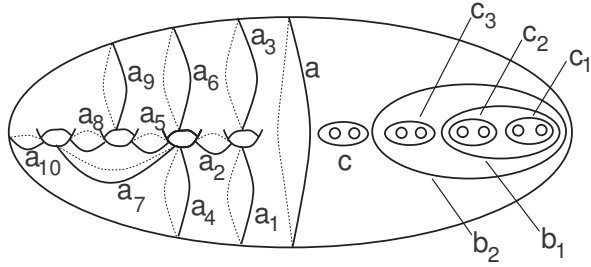} \caption{A $(0,3)$ circle,
$c$, in a pants decomposition}
\end{center}
\end{figure}

Let $a$ be a $(g,1)$-separating circle on $R$ which is disjoint
from $c$. Let $R_1$ be a subsurface of $R$ of genus $g$ having $a$
as its boundary and $R_2 =  R \setminus R_1 \cup a$. Let's choose
a pair of pants decomposition $P_1$ of $R_1$ as in Lemma
\ref{division}. If the $g =2$ we choose $P_1$ as in Figure 1
(i), if $g \geq 3$, we choose $P_1$ as in Figure 1, (ii).\\

Assume that $p=2n$ for some $n \in \mathbb{N}$ and $n \geq 2$.\\

(i) If $n = 2$, then we choose a pants decomposition $P$ on $R$
s.t. $P = P_1 \cup \{a, c, c_1 \}$ where $c_1$ is a $(0,3)$ curve
in $P$. Let $P'$ be a pair of pants decomposition of $R$ such that
$\lambda([P]) = [P']$. Let $a',c', c_1'$ be the representatives of
$\lambda([a]), \lambda([c]), \lambda([c_1])$ in $P'$ respectively.
Let $P_2'= \{c', c_1'\}$. By Lemma \ref{division} and Lemma 2.3,
there exist subsurfaces $R_1', R_2'$ of $R$ of genus $g$ and $0$
respectively s.t. $R_1'$ has $a'$ as its boundary, $R_2' = R
\setminus R_1' \cup a'$ and $P_2'$ is a pants decomposition for
$R_2'$. By using Lemma 2.4, we can see that $R_1' \cup C(a')$ is a
genus $g$ surface having $c', c_1'$ as its boundary components.
Since $R_1' \cup C(a')$ contains $P'$ and $P'$ is a pants
decomposition of $R$, each of $c', c_1'$
has to be a $(0,3)$ curve.\\

(ii) If $n > 2$, then we choose a pants decomposition $P$ on $R$
s.t. $P=P_1 \cup \{a, c, b_1, ..., b_{n-2}, c_1, ..., c_{n-1} \}$
where, $b_1, ..., b_{n-2}$ are 4-curves and $c_1, ..., c_{n-1}$
are $(0,3)$ curves in $P$, as shown in Figure 2 (for $g=4, p=8$).
In the other cases, a similar pair of pants decomposition of $R$
can be chosen. Let $P'$ be a pair of pants decomposition of $R$
such that $\lambda([P]) = [P']$. Let $a_i'$ be the representative
of $\lambda([a_i])$ in $P'$ for $i=1,..., 3g-2$, $a',c'$ be the
representatives of $\lambda([a]), \lambda([c])$ in $P'$
respectively, $b_i'$ be the representative of $\lambda([b_i])$ in
$P'$ for $i=1,..., n-2$, $c_i'$ be the representative of
$\lambda([c_i])$ in $P'$ for $i=1,..., n-1$. Let $P_1'= \{a_1',
a_2', ..., a_{3g-2}'\}$, $P_2' = P' \setminus (P_1' \cup
\{a'\})$.\\

By Lemma \ref{division} and Lemma 2.3, there exist subsurfaces
$R_1', R_2'$ of $R$ of genus $g$ and $0$ respectively s.t. $R_1'$
has $a'$ as its boundary, $R_2' = R \setminus R_1' \cup a'$ and
$P_1', P_2'$ are pants decompositions for $R_1', R_2'$
respectively. By using Lemma 2.4, we can see that $R_1' \cup
C(b'_1) \cup ... \cup C(b'_{n-2})$ is a genus $g$ surface having
$c', c_1', ..., c_{n-1}'$ as its boundary components. Since $R_1'
\cup C(b'_1) \cup ... \cup C(b'_{n-2})$ contains $P'$ and $P'$ is
a pants decomposition of $R$, each of $c', c_1', ..., c_n'$ has to
be a $(0,3)$ curve. Hence, $c'$ is a $(0,3)$ curve.\\

Assume that $p=2n+1$ for some $n \in \mathbb{N}$ and $n \geq 1$.\\

(i) If $n = 1$, then $P = P_1 \cup \{a, c \}$ is a pants
decomposition on $R$. Let $P'$ be a pair of pants decomposition of
$R$ such that $\lambda([P]) = [P']$. Let $a',c'$ be the
representatives of $\lambda([a]), \lambda([c])$ in $P'$
respectively. Let $P_2'= \{c'\}$. By Lemma \ref{division} and
Lemma 2.3, there exist subsurfaces $R_1', R_2'$ of $R$ of genus
$g$ and $0$ respectively s.t. $R_1'$ has $a'$ as its boundary,
$R_2' = R \setminus R_1' \cup a'$ and $P_2'$ is a pants
decomposition for $R_2'$. Since $a'$ is adjacent to $c'$ in $P'$,
there exists a pair of pants $Q$ in $R_2'$ having $a'$ and $c'$ on
its boundary. By using Lemma 2.4, we can see that $R_1' \cup Q$ is
a genus $g$ surface with two boundary components. One of the
boundary components is $c'$. Since $R_1' \cup Q$ contains $P'$ and
$P'$ is a pants decomposition of $R$, $c'$
has to be a $(0,3)$ curve.\\

\begin{figure}[htb]
\label{????????}
\begin{center}
\epsfxsize=3.33in \epsfbox{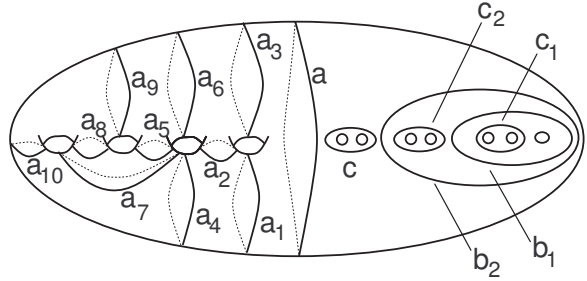} \caption{A $(0,3)$ circle,
$c$, in a pants decomposition}
\end{center}
\end{figure}

(ii) If $n = 2$, then we choose a pants decomposition $P$ on $R$
s.t. $P = P_1 \cup \{a, c, c_1, b_1 \}$ where $b_1, c_1$ are as in
Figure 3. Let $P'$ be a pair of pants decomposition of $R$ such
that $\lambda([P]) = [P']$. Let $a',c', c_1', b_1'$ be the
representatives of $\lambda([a]), \lambda([c]), \lambda([c_1]),
\lambda([b_1])$ in $P'$ respectively. Let $P_2'= \{c', c_1',
b_1'\}$. By Lemma \ref{division} and Lemma 2.3, there exist
subsurfaces $R_1', R_2'$ of $R$ of genus $g$ and $0$ respectively
s.t. $R_1'$ has $a'$ as its boundary, $R_2' = R \setminus R_1'
\cup a'$ and $P_2'$ is a pants decomposition for $R_2'$. By using
Lemma 2.4, we can see that $R_1' \cup C(a')$ is a genus $g$
surface having $c', b_1'$ as its boundary components. Since $b_1$
is adjacent to $c_1$ in $P$, $b_1'$ is adjacent to $c_1'$ in $P'$,
then there exists a pair of pants $Q$ having $b_1'$ and $c_1'$ on
its boundary. Then, since $R_1' \cup C(a') \cup Q$ contains $P'$
and $P'$ is a pants decomposition of $R$, each of $c', c_1'$
has to be a $(0,3)$ curve.\\

(iii) If $n > 2$, then we choose a pair of pants decomposition $P$
on $R$ s.t. $P=P_1 \cup \{a , c, b_1, ..., b_{n-1}, c_1, ...$,
$c_{n-1} \}$ where $P_1$ is a pair of pants decomposition as
before, $b_2, ..., b_{n-1}$ are 4-curves, $b_1$ is a 3-curve and
$c_1, ..., c_{n-1}$ are $(0,3)$ curves in $P$ as shown in Figure 3
(for $g=4, p=7$). Let $P'$ be a pair of pants decomposition of $R$
such that $\lambda([P]) = [P']$. Let $a_i'$ be the representative
of $\lambda([a_i])$ in $P'$ for $i=1,..., 3g-2$, $a',c'$ be the
representatives of $\lambda([a]), \lambda([c])$ in $P'$
respectively, $b_i'$ be the representative of $\lambda([b_i])$ in
$P'$ for $i=1,..., n-1$, $c_i'$ be the representative of
$\lambda([c_i])$ in $P'$ for $i=1,..., n-1$. Let $P_1'= \{a_1',
a_2', ..., a_{3g-2}'\}$, $P_2' = P' \setminus (P_1' \cup \{a\})$.
Let $R_1', R_2'$ be the subsurfaces of $R$ of genus $g$ and $0$
respectively s.t. $R_1'$ has $a'$ as its boundary, $R_2' = R
\setminus R_1' \cup a'$ and $P_2'$ is a pants decomposition for
$R_2'$. By using Lemma 2.4, we can see that $R_1' \cup C(b'_2)
\cup ... \cup C(b'_{n-1})$ is a genus $g$ subsurface having $c',
c_2', c_3',..., c_{n-1}', b_1'$ as its boundary components. Notice
that $b_1'$ is adjacent to $b_2', c_2'$ in $P'$ and $b_2', c_2'$
live in $R_1' \cup C(b'_2) \cup ... \cup C(b'_{n-1})$. Then, since
$b_1$ is adjacent to $c_1$ w.r.t. $P$, there is a pair of pants
$Q$ having $b_1'$ and $c_1'$ on its boundary. Then, $R_1' \cup
C(b'_2) \cup ... \cup C(b'_{n-1}) \cup Q$ is a genus $g$
subsurface having $c', c_1', c_2',..., c_{n-1}'$ as its boundary
components. Since $R_1' \cup C(b'_2) \cup ... \cup C(b'_{n-1})
\cup Q$ contains $P'$ and $P'$ is a pants decomposition of $R$,
each of $c', c_1', ..., c_{n-1}'$ has to be a $(0,3)$ curve.
Hence, $c'$ is a $(0,3)$ curve. \end{proof}\\

Let $\alpha$, $\beta$ be two distinct vertices in
$\mathcal{C}(R)$. We call $(\alpha, \beta)$ to be a peripheral
pair in $\mathcal{C}(R)$ if they have disjoint nonseparating
representatives $x, y$ respectively such that $x, y$ and a
boundary component of $R$ bound a pair of pants in $R$.

\begin{lemma}
\label{??} Let $\lambda : \mathcal{C}(R) \rightarrow
\mathcal{C}(R)$ be a superinjective simplicial map and $(\alpha,
\beta)$ be a peripheral pair in $\mathcal{C}(R)$. Then,
$(\lambda(\alpha), \lambda(\beta))$ is a peripheral pair in
$\mathcal{C}(R)$.
\end{lemma}

\begin{proof} Let $x, y$ be disjoint nonseparating representatives
of $\alpha, \beta$ respectively such that $x, y$ and a boundary
component of $R$ bound a pair of pants in $R$.\\

We will first prove the lemma when $g \geq 3$. Let $a$ be a
$(g-1,1)$-separating circle which is disjoint from $x$ and $y$.
Let $R_1$ be a subsurface of $R$ of genus $g-1$ having $a$ as its
boundary and $R_2 =  R \setminus R_1 \cup a$. Let's choose a pair
of pants decomposition $P_1$ of $R_1$ as in Lemma 2.5 (notice that
$P_1$ was chosen depending on $g$ in Lemma 2.5).\\

\begin{figure}[htb]
\label{?????}
\begin{center}
\epsfxsize=3.45in \epsfbox{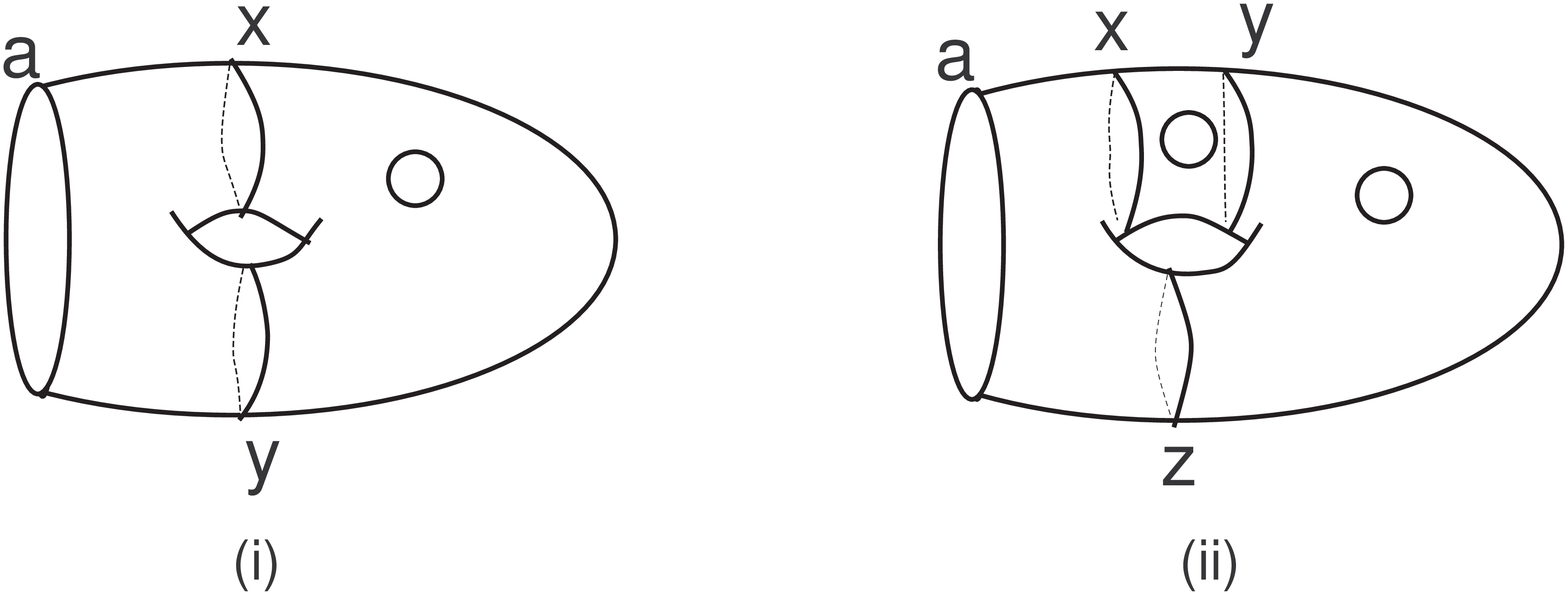} \caption{Peripheral pairs}
\end{center}
\end{figure}

Assume that $p=1$. Then $P = P_1 \cup \{a, x, y \}$ is a pair of
pants decomposition on $R$ (see Figure 4, (i)). Let $P'$ be a pair
of pants decomposition of $R$ such that $\lambda([P]) = [P']$. Let
$a', x', y'$ be the representatives of $\lambda([a]),
\lambda([x]), \lambda([y])$ in $P'$ respectively. Let $P_1'$ be
the set of elements in $P'$ which corresponds to $P_1$, and let
$P_2' = P' \setminus (P_1' \cup \{a'\})$. Since $g \geq 3$, by
Lemma 2.5, there exist subsurfaces $R_1', R_2'$ of $R$ of genus
$g-1$ and $1$ respectively s.t. $R_1'$ has $a'$ as its boundary,
$R_2' = R \setminus R_1' \cup a'$ and $P_2'= \{x', y'\}$ is a
pants decomposition for $R_2'$. Since $a$ is adjacent to $x$ and
$y$ w.r.t. $P$ in $R_2$, $a'$ has to be adjacent to $x'$ and $y'$
w.r.t. $P'$ in $R_2'$. So, there exists a pair of pants $Q$ in
$R_2'$ which has $a', x', y'$ on its boundary. Then, $R_1' \cup Q$
is a genus $g-1$ surface with two boundary components $x', y'$.
Then, since $R_1' \cup Q$ contains $P'$ and $P'$ is a pants
decomposition of $R$, and $R$ is a genus $g$ surface with 1
boundary component, there has to be a pair of pants having $x',
y'$ and the boundary component of $R$ on its
boundary. This proves the lemma when $p=1$.\\

Assume that $p=2$. Then $P = P_1 \cup \{a, x, y, z\}$ (see Figure
4, (ii)) is a pair of pants decomposition on $R$ where $z$ is a
nonseparating curve, $a, x, z$ are the boundary components of a
pair of pants, and $([y], [z])$ is a peripheral pair. Let $P'$ be
a pair of pants decomposition of $R$ such that $\lambda([P]) =
[P']$. Let $a', x', y', z'$ be the representatives of
$\lambda([a]), \lambda([x]), \lambda([y]), \lambda([z])$ in $P'$
respectively. Let $P_1'$ be the set of elements in $P'$ which
corresponds to $P_1$, and let $P_2' = P' \setminus (P_1' \cup
\{a'\})$.\\

By Lemma 2.5, there exist subsurfaces $R_1', R_2'$ of $R$ of genus
$g-1$ and $1$ respectively s.t. $R_1'$ has $a'$ as its boundary,
$R_2' = R \setminus R_1' \cup a'$ and $P_2'= \{x', y', z'\}$ is a
pants decomposition for $R_2'$. Since $a$ is adjacent to $x$ and
$z$ w.r.t. $P$ in $R_2$, $a'$ has to be adjacent to $x'$ and $z'$
w.r.t. $P'$ in $R_2'$. So, there exists a pair of pants $Q$ in
$R_2'$ which has $a', x', z'$ on its boundary. Then, $R_1' \cup Q$
is a genus $g-1$ surface with two boundary components $x', z'$.
Then, since $x$ is adjacent to $y$ in $R_2$, $x'$ is adjacent to
$y'$ in $R_2'$. So, there exists a pair of pants $T$ having $x',
y'$ on its boundary. Then, $R_1' \cup Q \cup T$ is a genus $g-1$
surface with three boundary components. Let $w$ be the boundary
component of $T$ which is different from $z', y'$. Since $R_1'
\cup Q \cup T$ contains $P'$ and $P'$ is a pants decomposition of
$R$, and $z'$ and $y'$ are nonseparating circles by Lemma 2.6, $w$
is a boundary component of $R$. This proves the lemma for
$p=2$.\\

\begin{figure}[htb]
\begin{center}
\epsfxsize=5in \epsfbox{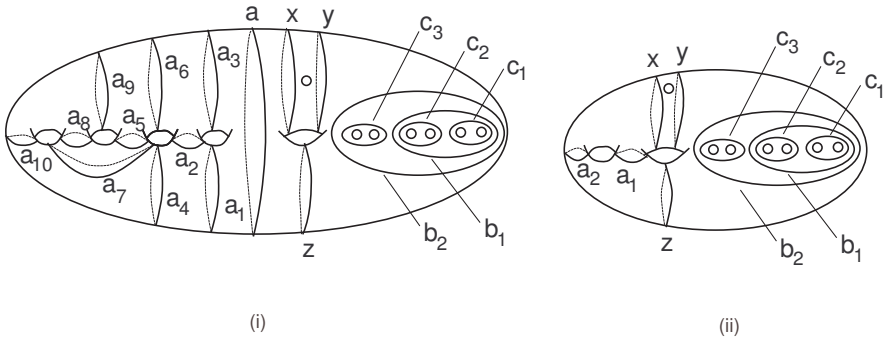} \caption{Peripheral pairs
pants decompositions}
\end{center}
\end{figure}

Assume that $p=2n+1$ and $n \geq 1$.\\

(i) If $n =1$, we choose a pants decomposition $P$ on $R$ s.t. $P
= P_1 \cup \{a, x, y, z, c_1 \}$ where $z$ is a 4-curve, $c_1$ is
a $(0,3)$ curve in $P$, $a, x, z$ bound a pair of pants in $P$ and
$z, y, c_1$ bound a pair of pants in $P$. Let $P'$ be a pair of
pants decomposition of $R$ such that $\lambda([P]) = [P']$. Let
$a', x', y', z', c_1'$ be the representatives of $\lambda([a]),
\lambda([x]), \lambda([y]), \lambda([z]), \lambda([c_1])$ in $P'$
respectively. Let $P_1'$ be the set of elements in $P'$ which
corresponds to $P_1$, and let $P_2' = \{x', y', z', c_1' \} $.\\

By Lemma 2.5, there exist subsurfaces $R_1', R_2'$ of $R$ of genus
$g-1$ and $1$ respectively s.t. $R_1'$ has $a'$ as its boundary,
$R_2' = R \setminus R_1' \cup a'$ and $P_2'= \{x', y', z', c_1'\}$
is a pants decomposition for $R_2'$. Since $a$ is adjacent to $x$
and $z$ w.r.t. $P$ in $R_2$, $a'$ has to be adjacent to $x'$ and
$z'$ w.r.t. $P'$ in $R_2'$. So, there exists a pair of pants $Q$
in $R_2'$ which has $a', x', z'$ on its boundary. Then, $R_1' \cup
Q$ is a genus $g-1$ surface with two boundary components $x', z'$.
Then, since $z$ is adjacent to $y$ and $c_1$ w.r.t. $P$, $z'$ is
adjacent to $y'$ and $c_1'$ w.r.t. $P'$ in $R_2'$. So, there
exists a pair of pants $T$ in $R_2'$ which has $z', y', c_1'$ on
its boundary. Then, $R_1' \cup Q \cup T$ is a genus $g-1$ surface
with three boundary components $x', y', c_1'$. Since $c_1$ is a
$(0,3)$ curve in $P$, $c_1'$ is a $(0,3)$ curve by Lemma 2.7.
Then, since $R_1' \cup Q \cup T$ contains $P'$ and $P'$ is a pants
decomposition of $R$, there has to be a pair of pants containing
$x'$, $y'$ and a boundary component of $R$ as its boundary components.\\

(ii) If $n > 1$, we choose a pair of pants decomposition $P$ on
$R$ such that $P = P_1 \cup \{a, x, y, z, b_1,..., b_{n-1}, c_1,
..., c_n \}$ where $z, b_1, ..., b_{n-1}$ are 4-curves and $c_1,
..., c_n$ are $(0,3)$ curves in $P$ and $a, x, z$ bound a pair of
pants in $P$. In Figure 5 (i), we show how to choose $P$ when
$g=5, p=7$. In the other cases, a similar pair of pants
decomposition of $R$ can be chosen. Let $P'$ be a pair of pants
decomposition of $R$ such that $\lambda([P]) = [P']$. Let $a', x',
y', z'$ be the representatives of $\lambda([a]), \lambda([x]),
\lambda([y]), \lambda([z])$ in $P'$ respectively, $b_i'$ be the
representative of $\lambda([b_i])$ in $P'$ for $i=1,..., n-1$ and
$c_i'$ be the representative of $\lambda([c_i])$ in $P'$ for
$i=1,..., n$. Let $P_1'$ be the set of elements in $P'$ which
corresponds to $P_1$, and let $P_2' = P' \setminus
(P_1' \cup \{a'\})$.\\

By Lemma \ref{division} and Lemma 2.3, there exist subsurfaces
$R_1', R_2'$ of $R$ of genus $g-1$ and $1$ respectively s.t.
$R_1'$ has $a'$ as its boundary, $R_2' = R \setminus R_1' \cup a'$
and $P_2'$ is a pants decomposition for $R_2'$. We can see that
$R_1' \cup C(z') \cup C(b'_1) \cup ... \cup C(b'_{n-1})$ is a
genus $g-1$ surface having $x', y', c_1', ..., c_n'$ as its
boundary components. Since $c_1, ..., c_n$ are $(0,3)$ curves in
$P$, $c_1', ..., c_n'$ are $(0,3)$ curves by Lemma 2.7. Then,
since $R_1' \cup C(z') \cup C(b'_1) \cup ... \cup C(b'_{n-1})$
contains $P'$ and $P'$ is a pants decomposition of $R$, there has
to be a pair of pants containing $x'$ and $y'$ and a boundary
component of $R$ as its boundary components.\\

\begin{figure}[htb]
\label{????}
\begin{center}
\epsfxsize=5in \epsfbox{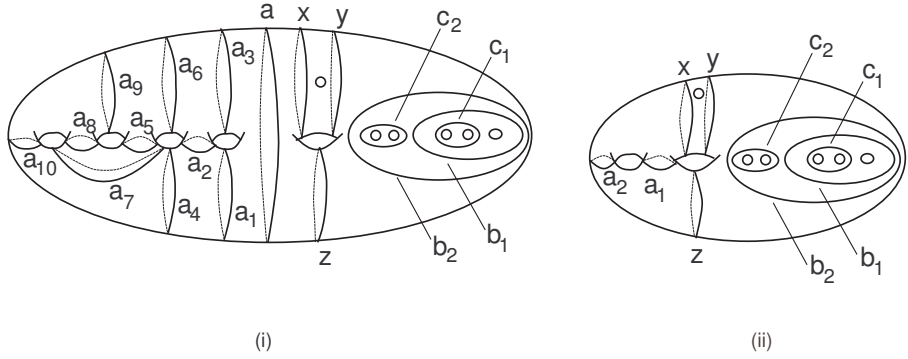} \caption{Peripheral pairs in
pants decompositions}
\end{center}
\end{figure}

Assume that $p=2n$, and $n \geq 2$. Let's choose a pants
decomposition $P$ on $R$ s.t. $P = P_1 \cup \{a, x, y, z, b_1,
..., b_{n-1}, c_1, ..., c_{n-1} \}$ where, $z, b_2, ..., b_{n-1}$
are 4-curves, $b_1$ is a 3-curve and $c_1, ..., c_n$ are $(0,3)$
curves in $P$. In Figure 6 (i), we show how to choose $P$ when
$g=5, p=6$. Let $P'$ be a pair of pants decomposition of $R$ such
that $\lambda([P]) = [P']$. Let $a', x', y', z'$ be the
representatives of $\lambda([a]), \lambda([x]), \lambda([y]),
\lambda([z])$ in $P'$ respectively, $b_i'$ be the representative
of $\lambda([b_i])$ in $P'$ for $i=1,..., n-1$, $c_i'$ be the
representative of $\lambda([c_i])$ in $P'$ for $i=1,..., n-1$. Let
$P_1'$ be the set of elements in $P'$ which corresponds to $P_1$,
and let $P_2' = P' \setminus (P_1' \cup \{a'\})$. By Lemma
\ref{division} and Lemma 2.3, there exist subsurfaces $R_1', R_2'$
of $R$ of genus $g-1$ and $1$ respectively s.t. $R_1'$ has $a'$ as
its boundary, $R_2' = R \setminus R_1' \cup a'$
and $P_2'$ is a pants decomposition for $R_2'$.\\

(i) If $n=2$, $V = R_1' \cup C(z')$ is a genus $g-1$ surface
having $x', y', b_1'$ as its boundary components. Since $b_1$ is
adjacent to $c_1$ w.r.t. $P$, $b_1'$ is adjacent to $c_1'$ w.r.t.
$P'$, so there exists a pair of pants $Y$ in $R \setminus V$ s.t.
$Y$ contains $b_1', c_1'$ on its boundary. Since $c_1$ is a
$(0,3)$ curve in $P$, $c_1'$ is a $(0,3)$ curve by Lemma 2.7.
Then, since $V \cup Y$ contains $P'$ and $P'$ is a pants
decomposition of $R$, and $x'$ and $y'$ are nonseparating circles,
there has to be a pair of pants containing $x'$ and $y'$ and a
boundary component of
$R$ as its boundary components.\\

(ii) If $n > 2$, we can see that $W = R_1' \cup C(z') \cup C(b'_2)
\cup ... \cup C(b'_{n-1})$ is a genus $g-1$ surface having $x',
y', c_2', ..., c_{n-1}', b_1'$ as its boundary components. Since
$b_1$ is adjacent to $c_1$ w.r.t. $P$, $b_1'$ is adjacent to
$c_1'$ w.r.t. $P'$, so there exists a pair of pants $Y$ in $R
\setminus W$ s.t. $Y$ contains $b_1', c_1'$ on its boundary. Since
$c_1, ..., c_{n-1}$ are $(0,3)$ curves in $P$, $c_1', ...,
c_{n-1}'$ are $(0,3)$ curves by Lemma 2.7. Then, since $W \cup Y$
contains $P'$ and $P'$ is a pants decomposition of $R$, and $x'$
and $y'$ are nonseparating circles, there has to be a pair of
pants containing $x'$ and $y'$ and a boundary component of $R$ as
its boundary
components. This proves the lemma when $g \geq 3$.\\

When $g = 2$ the proof is similar. We have $p \geq 2$ in this
case. Instead of using a separating curve $a$, we use pair of
pants decompositions as given (for special cases) in Figure 5 (ii)
and Figure 6 (ii). By the proof of Lemma 2.5, there are pairwise
disjoint representatives $a_1', a_2', z', x'$ of $\lambda([a_1]),
\lambda([a_2]), \lambda([z]), \lambda([x])$ and a subsurface
$R_1'$ of genus 1 with two boundary components $x'$ and $z'$ such
that $x'$, $a_1'$, $a_2'$ bound a pair of pants in $R_1'$, and
$z'$, $a_1'$, $a_2'$ bound a pair of pants in $R_1'$. Then, we
follow the proof of the case when $g \geq 3$.\end{proof}\\

Let $M$ be a sphere with $k$ holes and $k \geq 5$. A circle $a$ on
$M$ is called an {\it n-circle} if $a$ bounds a disk with $n$
holes on $M$ where $n \geq 2$. If $a$ is a 2-circle on $M$, then
there exists up to isotopy a unique nontrivial embedded arc $a'$
on the two-holed disk component of $M_a$ joining the two holes in
this disk. If $a$ and $b$ are two 2-circles on $M$ such that the
corresponding arcs $a'$, $b'$ can be chosen to meet exactly at one
common end point, and $\alpha = [a], \beta = [b]$, then $\{
\alpha, \beta \}$ is called a {\it simple pair}. A {\it pentagon}
in $\mathcal{C}(M)$ is an ordered 5-tuple $(\alpha_1, \alpha_2,
\alpha_3, \alpha_4, \alpha_5)$, defined up to cyclic permutations,
of vertices of $\mathcal{C}(M)$ such that $i(\alpha_j,
\alpha_{j+1}) = 0$ for $j=1,2,...,5$ and $i(\alpha_j, \alpha_k)
\neq 0$ otherwise, where $\alpha_6 = \alpha_1$. A vertex in
$\mathcal{C}(M)$ is called an {\it n-vertex} if it has a
representative which is an n-circle on $M$. Let $M'$ be the
interior of $M$. There is a natural isomorphism $\chi:
\mathcal{C}(M') \rightarrow \mathcal{C}(M)$ which respects the
above notions and the corresponding notions in \cite{K}. Using
this isomorphism, we can restate a theorem of Korkmaz as follows:

\begin{theorem}
\label{korkmaz} {\bf (Korkmaz) \cite{K}} Let $M$ be a sphere with
$n$ holes and $n \geq 5$. Let $\alpha, \beta$ be two 2-vertices of
$\mathcal{C(M)}$. Then $\{ \alpha, \beta \}$ is a simple pair iff
there exist vertices $\gamma_1, \gamma_2, ..., \gamma_{n-2}$ of
$\mathcal{C}(M)$ satisfying the following conditions: \\
\indent (i) $(\gamma_1, \gamma_2, \alpha, \gamma_3, \beta)$ is a
pentagon
in $\mathcal{C}(M)$, \\
\indent (ii) $\gamma_1$ and $\gamma_{n-2}$ are 2-vertices,
$\gamma_2$ is a 3-vertex and $\gamma_k$ and $\gamma_{n-k}$ are
k-vertices for $3
\leq k \leq \frac{n}{2}$,\\
\indent (iii) $\{\alpha, \gamma_3, \gamma_4, \gamma_5, ...,
\gamma_{n-2}\}$, $\{\alpha, \gamma_2, \gamma_4, \gamma_5, ...,
\gamma_{n-2}\}$,  $\{\beta, \gamma_3, \gamma_4, \gamma_5, ...,
\gamma_{n-2}\}$, and $\{\gamma_1, \gamma_2, \gamma_4, \gamma_5,
..., \gamma_{n-2}\}$ are codimension-zero simplices.\end{theorem}

\begin{lemma}
\label{separating} Let $\lambda : \mathcal{C}(R) \rightarrow
\mathcal{C}(R)$ be a superinjective simplicial map. Then,
$\lambda$ sends the isotopy class of a $(k,m)$-separating circle
to the isotopy class of a $(k,m)$-separating circle, where $1 \leq
k \leq g$, $ 1 \leq m < p$.\end{lemma}

\begin{proof} Let $\alpha = [a]$ where $a$ is a $(k,m)$-separating
circle where $1 \leq k \leq g$, $ 1 \leq m < p$.\\

Case 1: Assume the genus of $R$ is at least 3. Then $a$ separates
a subsurface of genus at least 2. So, it is enough to consider the
cases when $k \geq 2$. If $m=1$, then the lemma follows by Lemma
2.5. Assume that $m \geq 2$. Let $R_1$ be a subsurface of $R$ of
genus $k$ with $m$ boundary components which has $a$ as one of its
boundary components. Let's choose a pair of pants decomposition
$P_1= \{a_1, a_2,... a_{3k-3}, b_1,...,b_m\}$ of $R_1$ where
$a_1,.., a_{3k-3}$ are 4-curves and $(b_i, b_{i+1})$ is a
peripheral pair for $i=1,...,m-1$ as shown in Figure 7 (i) (for
$k=3, m=5$). Then, we complete $P_1 \cup \{a\}$ to a pair of pants
decomposition $P$ of $R$ in any way we like. By Lemma
\ref{imageofpantsdecomp}, we can choose a pair of pants
decomposition, $P'$, of $R$ such that $\lambda([P]) = [P']$.\\

\begin{figure}[htb]
\begin{center}
\epsfxsize=6in \epsfbox{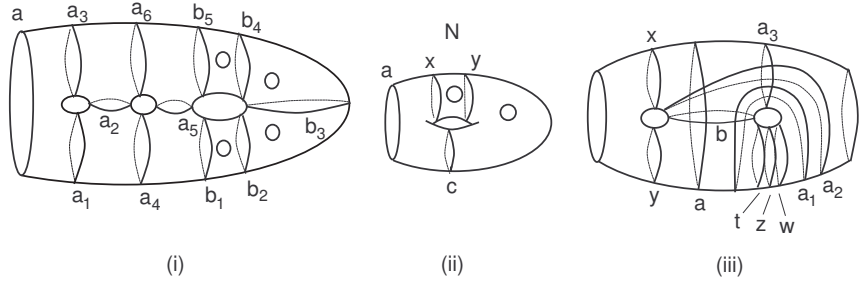} \caption{(i) A $(3,5)$
circle $a$, (ii) A $(1,3)$ circle $a$, (iii) A $(1,2)$ circle $a$}
\end{center}
\end{figure}

Let $a_i'$ be the representative of $\lambda([a_i])$ which is in
$P'$ for $i=1,..,3k-3$ and $a'$ be the representative of
$\lambda([a])$ which is in $P'$. Let $b_i'$ be the representative
of $\lambda([b_i])$ which is in $P'$ for $i=1,..,m$. Since $a_i$
is a 4-curve in $P$, for $i=1,..,3k-3$, by Lemma \ref{adjacent}
and Lemma \ref{injective}, $a_i'$ is a 4-curve in $P'$, for
$i=1,..,3k-3$. Since $(b_i, b_{i+1})$ is a peripheral pair,
$(b_i', b_{i+1}')$ is a peripheral pair for $i=1,...,m-1$ by Lemma
2.8. Then, there exist distinct pair of pants $Q_i$ which has
$b_i, b_{i+1}$ and a boundary component of $R$ on its boundary for
$i=1,...,m-1$. Then, it is easy to see that $C(a_1') \cup ...
C(a_{3k-3}') \cup Q_1 \cup ... \cup Q_{m-1}$ is a genus $k$
subsurface with $m$ boundary components having $a'$ on the
boundary. All the other boundary components of this subsurface are
boundary components of $R$. Hence, $a'$ is a $(k,
m)$ separating curve.\\

Case 2: Assume the genus of $R$ is 2 and the number of boundary
components is at least 2. If $a$ separates a subsurface of genus
2, then the proof is similar to the proof of case 1. Assume that
$a$ separates a subsurface of genus 1. Then, since the number of
boundary components is at least 2, $a$ separates a subsurface of
genus 1 with at least 2 boundary components. We will consider the
following two cases:\\

(i) $a$ separates a subsurface $N$ of genus 1 with at least 3
boundary components. We will give the proof when $a$ is a $(1,3)$
circle. The proofs of the remaining cases are similar. Let $N$ and
$x, y, c$ be as shown in Figure 7 (ii). We complete $\{a, x, y,
c\}$ to a pair of pants decomposition $P$ of $R$. Let $P'$ be a
pair of pants decomposition of $R$ such that $\lambda([P]) =
[P']$. Let $a', x', y', c'$ be the representatives of
$\lambda([a]), \lambda([x]), \lambda([y]), \lambda([c])$ in $P'$
respectively. Since $(x, y)$ and $(c, y)$ are peripheral pairs,
$(x', y')$ and $(c', y')$ are peripheral pairs by Lemma 2.8. Then,
there exist distinct pair of pants $Q_1, Q_2$ such that $Q_1$ has
$x', y'$ and a boundary component of $R$ on its boundary, and
$Q_2$ has $c', y'$ and a boundary component of $R$ on its
boundary. There exists also a pair of pants $Q_3$ which has $x',
a', c'$ on its boundary since $x'$ is adjacent to $a'$ and $c'$
w.r.t. $P'$. Then, it is easy to see that $Q_1 \cup Q_2 \cup Q_3$
is a genus 1 subsurface with 3 boundary components having $a'$ as
one of its boundary components and all the other boundary
components of this subsurface are boundary components of $R$.
Hence, $a'$ is a $(1, 3)$ separating curve.\\

(ii) Each connected component of $R_a$ is a genus 1 surface with 2
boundary components. Let's choose a pair of pants decomposition $P
= \{x, y, z, a_3, b\}$ on $R$, where the curves are as in Figure 7
(iii). Note that $b$ and $a$ have geometric intersection 2,
algebraic 0. Let $P'$ be a pair of pants decomposition of $R$ such
that $\lambda([P]) = [P']$. Let $x', y', z', a_3', b'$ be the
representatives of $\lambda([x]), \lambda([y]), \lambda([z]),
\lambda([a_3]), \lambda([b])$ in $P'$ respectively. Since $(x, y)$
and $(a_3, z)$ are peripheral pairs, $(x', y')$ and $(a_3', z')$
are peripheral pairs by Lemma 2.8. Then, there exist distinct pair
of pants $Q_1', Q_2'$ such that $Q_1'$ has $x', y'$ and a boundary
component of $R$ on its boundary, and $Q_2'$ has $a_3', z'$ and a
boundary component of $R$ on its boundary. Then, since $x'$ is
adjacent to $a_3'$ and $b'$ w.r.t. $P'$, and $y'$ is adjacent to
$z'$ and $b'$ w.r.t. $P'$, there exist also distinct pair of pants
$Q_3'$ and $Q_4'$ having $x', a_3', b'$ and $y', z', b'$ on their
boundary respectively. Then, it is easy to see that $R = Q_1' \cup
Q_2' \cup Q_3' \cup Q_4'$ and there exists a homeomorphism $\chi :
(R, x, y, b, z, a_3) \rightarrow (R, x', y', b', z', a_3')$, i.e.
$P$ and $P'$ are topologically equivalent.\\


Let $a_1, a_2, t, w$ be as shown in Figure 7 (iii). Let $Q_2$ be
the pair of pants in $P$, with boundary components $z, a_3$ and a
boundary component of $R$. Let $M$ be the subsurface of $R$
bounded by $x, y, t, a_3$. $M$ is a sphere with four holes. Let
$\tilde{M}$ be the subsurface of $R$ bounded by $x, y, t, w$ and
the boundary component of $R$ which is
on the boundary of $Q_2$. $\tilde{M}$ is a sphere with five holes.\\

Let $t' = \chi(t)$ and $w'= \chi(w)$. Let $M'$ be the subsurface
of $R$ bounded by $x', y', t', a_3'$. $M'$ is a sphere with four
holes. Let $\tilde{M'}$ be the subsurface of $R$ bounded by $x',
y', t', w'$ and the boundary component of $R$ which is on the
boundary of $Q_2'$. $\tilde{M'}$ is a sphere with
five holes.\\

Since $x, y, t, w$ are essential circles in $R$, the essential
circles on $\tilde{M}$ are essential in $R$. Similarly, since $x',
y', t', w'$ are essential circles in $R$, the essential circles on
$\tilde{M'}$ are essential in $R$. Furthermore, we can identify
$\mathcal{C}(\tilde{M})$ and $\mathcal{C}(\tilde{M'})$ with two
subcomplexes of $\mathcal{C}(R)$ in such a way that the isotopy
class of an essential circle in $\tilde{M}$ or in $\tilde{M'}$ is
identified with the isotopy class of that circle in $R$. Now,
suppose that $\alpha$ is a vertex in $\mathcal{C}(\tilde{M})$.
Then, with this identification, $\alpha$ is a vertex in
$\mathcal{C}(R)$ and $\alpha$ has a representative in $\tilde{M}$.
Then, $i(\alpha, [x])=i(\alpha, [t])=i(\alpha, [w])= i(\alpha,
[y])=0$. Then there are two possibilities: (i) $\alpha = [b]$ or
$\alpha = [a_3]$, (ii) $i(\alpha, [b]) \neq 0$ or $i(\alpha,
[a_3]) \neq 0$. Since $\lambda$ is injective, $\lambda(\alpha)$ is
not equal to any of $[x'], [t'], [w'], [y']$. Since $\lambda$ is
superinjective, $i(\lambda(\alpha), [x'])=i(\lambda(\alpha),
[t'])=i(\lambda(\alpha), [w'])= i(\lambda(\alpha), [y'])=0$. Then,
there are two possibilities: (i) $\lambda(\alpha) = [b']$ or
$\lambda(\alpha) = [a_3']$, (ii) $i(\lambda(\alpha), [b']) \neq 0$
or $i(\lambda(\alpha), [a_3']) \neq 0$. Then, a representative of
$\lambda(\alpha)$ can be chosen in $\tilde{M'}$. Hence, $\lambda$
maps the vertices of $\mathcal{C}(R)$ that have essential
representatives in $\tilde{M}$ to the vertices of $\mathcal{C}(R)$
that have essential representatives in $\tilde{M'}$, (i.e.
$\lambda$ maps $\mathcal{C}(\tilde{M}) \subseteq \mathcal{C}(R)$
to $\mathcal{C}(\tilde{M'}) \subseteq \mathcal{C}(R)$). Similarly,
$\lambda$ maps $\mathcal{C}(M) \subseteq \mathcal{C}(R)$ to
$\mathcal{C}(M') \subseteq \mathcal{C}(R)$.\\

Let $\gamma_1, \gamma_2, \gamma_3$ be the isotopy classes of $a_1,
a_2, a_3$ in $\tilde{M}$ respectively. It is easy to see that
$\{[b], [a]\}$ is a simple pair in $\tilde{M}$, $(\gamma_1,
\gamma_2, [b], \gamma_3, [a])$ is a pentagon in
$\mathcal{C}(\tilde{M})$, $\gamma_1$ and $\gamma_3$ are
2-vertices, $\gamma_2$ is a 3-vertex, and $\{[b], \gamma_3\}$,
$\{[b], \gamma_2\}$, $\{[a], \gamma_3\}$ and $\{\gamma_1, \gamma_2
\}$ are codimension-zero simplices of $\mathcal{C}(\tilde{M})$.\\

Since $\lambda$ is superinjective and $x', y', t', w'$ are
essential circles, we can see that $(\lambda(\gamma_1),
\lambda(\gamma_2)$, $\lambda([b]), \lambda(\gamma_3)$,
$\lambda([a]))$ is a pentagon in $\mathcal{C}(\tilde{M'})$. By the
proof of case 2 (i), we can see that $\lambda(\gamma_1)$ is a
2-vertex in $\mathcal{C}(\tilde{M'})$. Using $\chi$, it is easy to
see that $\lambda(\gamma_3)$ is a 2-vertex in
$\mathcal{C}(\tilde{M'})$. Since $(x, a_2)$ is a peripheral pair
in $\tilde{M}$, by using Lemma 2.8 and the existence of $\chi$ we
can see that $(x', a_2')$ is a peripheral pair in $\tilde{M'}$ and
so $\lambda(\gamma_2)$ is a 3-vertex in $\mathcal{C}(\tilde{M'})$.
Since $\lambda$ is an injective simplicial map $\{\lambda([b]),
\lambda(\gamma_3)\}$, $\{\lambda([b]), \lambda(\gamma_2)\}$,
$\{\lambda([a]), \lambda(\gamma_3)\}$ and $\{\lambda(\gamma_1),
\lambda(\gamma_2)\}$ are codimension-zero simplices of
$\mathcal{C}(\tilde{M'})$. Then, by Theorem \ref{korkmaz},
$\{\lambda([b]), \lambda([a])\}$ is a simple pair in $\tilde{M'}$.
Since $\lambda$ maps $\mathcal{C}(M)$ to $\mathcal{C}(M')$,
$\lambda([a])$ has a representative $a'$ in $M'$ such that
$i(\lambda([b]), \lambda([a]) = |b' \cap a'|$. Then, there exists
a homeomorphism $\theta : (R, x, y, b, a, a_3, z) \rightarrow (R,
x', y', b', a', a_3', z')$. Hence, $\lambda([a])$ has a
representative $a'$ which is a $(1,2)$ circle. This proves the
lemma for case 2, (ii).\end{proof}

\begin{lemma}
\label{division2} Let $\lambda : \mathcal{C}(R) \rightarrow
\mathcal{C}(R)$ be a superinjective simplicial map. Let $t$ be a
$(k,m)$-separating circle on $R$, where $1 \leq k \leq g$, $1 \leq
m < p$ separating $R$ into two subsurfaces $R_1, R_2$. Let $t' \in
\lambda([t])$. Then $t'$ is a $(k,m)$-separating circle, $t'$
separates $R$ into two subsurfaces $R_1', R_2'$ such that
$\lambda(\mathcal{C}(R_1)) \subseteq \mathcal{C}(R_1')$ and
$\lambda(\mathcal{C}(R_2)) \subseteq \mathcal{C}(R_2')$.
\end{lemma}

\begin{proof} Let $t$ be a $(k,m)$-separating circle where $1 \leq k \leq g$,
$1 \leq m < p$. Let $R_1, R_2$ be the distinct subsurfaces of $R$
of genus $k$ and $g-k$ respectively which come from the separation
by $t$. Let $t' \in \lambda([t])$. By Lemma \ref{separating}, $t'$
is a $(k,m)$-separating circle. As we showed in the proof of Lemma
\ref{separating}, there is a pair of pants decomposition $P_1$ of
$R_1$, and $P_1 \cup \{t\}$ can be completed to a pair of pants
decomposition $P$ of $R$ such that a set of curves, $P_1'$,
corresponding (via $\lambda$) to the curves in $P_1$, can be
chosen such that $P_1'$ is a pair of pants decomposition of a
subsurface that has $t'$ as a boundary component and all the other
boundary components of this subsurface are boundary components of
$R$. Let $R_1'$ be this subsurface. Let $R_2' = R \setminus R_1'
\cup t'$. A pairwise disjoint representative set, $P'$, of
$\lambda([P])$ containing $P_1' \cup \{t'\}$ can be chosen. Then,
by Lemma \ref{imageofpantsdecomp}, $P'$ is a pair of pants
decomposition of $R$. Let $P_2 = P \setminus (P_1 \cup {t})$ and
$P_2' = P' \setminus (P_1' \cup {t'})$. Then $P_2, P_2'$ are pair
of pants decompositions of $R_2, R_2'$ respectively as $P_1, P_1'$
are pair of pants decompositions of $R_1, R_1'$
respectively.\\

Now, let $\alpha$ be a vertex in $\mathcal{C}(R_1)$. Then, either
$\alpha \in [P_1]$ or $\alpha$ has a nonzero geometric
intersection with an element of $[P_1]$. In the first case,
clearly $\lambda(\alpha) \in \mathcal{C}(R_1')$ since elements of
$[P_1]$ correspond to elements of $[P_1'] \subseteq
\mathcal{C}(R_1')$. In the second case, since $\lambda$ preserves
zero and nonzero geometric intersection (since $\lambda$ is
superinjective) and $\alpha$ has zero geometric intersection with
the elements of $[P_2]$ and $[t]$, and nonzero intersection with
an element of $[P_1]$, $\lambda(\alpha)$ has zero geometric
intersection with elements of $[P_2']$ and $[t']$, and nonzero
intersection with an element of $[P_1']$. Then, $\lambda(\alpha)
\in \mathcal{C}(R_1')$. Hence, $\lambda(\mathcal{C}(R_1))
\subseteq \mathcal{C}(R_1')$. The proof of
$\lambda(\mathcal{C}(R_2)) \subseteq \mathcal{C}(R_2')$ is
similar.\end{proof}

\begin{lemma}
\label{top} Let $\lambda : \mathcal{C}(R) \rightarrow
\mathcal{C}(R)$ be a superinjective simplicial map. Then $\lambda$
preserves topological equivalence of ordered pairs of pants
decompositions of $R$, (i.e. for a given ordered pair of pants
decomposition $P=(c_1, c_2, ..., c_{3g-3+p})$ of $R$, and a
corresponding ordered pair of pants decomposition $P'=(c_1', c_2',
...,$ $ c_{3g-3+p}')$ of $R$, where $[c_i']= \lambda([c_i])$
$\forall i= 1, 2, ..., 3g-3+p$, there exists a homeomorphism $H: R
\rightarrow R$ such that $H(c_i)=c_i'$ $\forall i= 1, 2, ...,
3g-3+p$).\end{lemma}

\begin{proof} Let $P$ be a pair of pants decomposition of $R$ and
$A$ be a nonembedded pair of pants in $P$. The boundary of $A$
consists of the circles $x, y$ where $x$ is a 1-separating circle
on $R$ and $y$ is a nonseparating circle on $R$. Let $R_1$ be the
subsurface of $R$ of genus $g-1$ with $p+1$ boundary components
which has $x$ as one of its boundary component and let $R_2$ be
the subsurface of $R$ of genus 1 which is bounded by $x$. Let
$P_1$ be the set of elements of $P \setminus \{x\}$ which are on
$R_1$ and $P_2$ be the set of elements of $P \setminus \{x\}$
which are on $R_2$. Then, $P_1, P_2$ are pair of pants
decompositions of $R_1, R_2$ respectively. So, $P_2 = \{y\}$ is a
pair of pants decomposition of $R_2$. By Lemma \ref{division2},
there exists a 1-separating circle $x' \in \lambda([x])$ and
subsurfaces $R_1', R_2'$, of $R$, of genus $g-1$ and 1
respectively such that $\lambda(\mathcal{C}(R_1)) \subseteq
\mathcal{C}(R_1')$ and $\lambda(\mathcal{C}(R_2)) \subseteq
\mathcal{C}(R_2')$. Since $[P_1] \subseteq \mathcal{C}(R_1)$, we
have $\lambda([P_1]) \subseteq \mathcal{C}(R_1')$. Since $[P_2]
\subseteq \mathcal{C}(R_2)$, we have $\lambda([P_2]) \subseteq
\mathcal{C}(R_2')$. Since $\lambda$ preserves disjointness, we can
see that a set, $P_1'$, of pairwise disjoint representatives of
$\lambda([P_1])$ disjoint from $x'$ can be chosen. By counting the
number of curves in $P_1'$, we can see that $P_1'$ is a pair of
pants decomposition of $R_1'$. Similarly, a set, $P_2'$, of
pairwise disjoint representatives of $\lambda([P_2])$ disjoint
from $x'$ can be chosen. By counting the number of curves in
$P_2'$,  we can see that $P_2'$ is a pair of pants decomposition
of $R_2'$. Since $P_2$ has one element, $y$, $P_2'$ has one
element. Let $y' \in P_2'$. Since $x', y'$ correspond to $x,y$
respectively and $y$ and $y'$ give pair of pants decompositions on
$R_2$ and $R_2'$ (which are both nonembedded pairs of pants) and
$x$ and $x'$ are the boundaries of $R_2$ and $R_2'$, we see that
$\lambda$ ``sends" a nonembedded pair of pants to a nonembedded
pair of pants.\\

Let $B$ be an embedded pair of pants of $P$. Let $x, y, z$ be the
boundary components of $B$. There are three cases to
consider:\\

(i) At least one of $x$, $y$ or $z$ is a separating circle.\\
\indent (ii) All of $x, y, z$ are nonseparating circles.\\
\indent (iii) Exactly one of $x, y, z$ is a boundary component of
$R$ and the other two are nonseparating circles.\\

Case (i): W.L.O.G assume that $x$ is a $(k,m)$-separating circle
for $1 \leq k \leq g$, $ 1 \leq m < p$. Let $R_1, R_2$ be the
distinct subsurfaces of $R$ of genus $k$ and $g-k$ respectively
which comes from separation by $x$. W.L.O.G. assume that $y, z$
are on $R_2$. Let $x' \in \lambda([x])$. By Lemma \ref{division2},
$x'$ is a $(k,m)$ circle separating $R$ into subsurfaces, $R_1',
R_2'$, of genus $k$ and $g-k$ respectively such that
$\lambda(\mathcal{C}(R_1)) \subseteq \mathcal{C}(R_1')$ and
$\lambda(\mathcal{C}(R_2)) \subseteq \mathcal{C}(R_2')$.\\

(a) If $y$ and $z$ are nontrivial circles, then, since $y \cup z
\subseteq R_2$, $\lambda(\{[y], [z]\}) \subseteq
\mathcal{C}(R_2')$. Let $y' \in \lambda([y]), z' \in \lambda([z])$
such that $\{x', y', z'\}$ is pairwise disjoint. Let $P'$ be a set
of pairwise disjoint representatives of $\lambda([P])$ which
contains $x', y', z'$. $P'$ is a pair of pants decomposition of
$R$. Then, since $x$ is adjacent to $y$ and $z$ w.r.t. $P$, $x'$
is adjacent to $y'$ and $z'$ w.r.t. $P'$ by Lemma \ref{adjacent}.
Then, since $x' \cup y'\cup z' \subseteq R_2'$, and $x'$ is a
boundary component of $R_2'$, there is an embedded pair of pants
in $R_2'$ which has $x', y', z'$ on its boundary.\\

(b) If each of $y$ and $z$ is a boundary component of $R$, then
$x$ is a (0,3)-separating circle. Then, by Lemma 2.7, there exists
$x' \in \lambda([x])$ s.t. $x'$ is a (0,3)-separating circle. So,
there is an embedded pair of pants which has $x'$ and two boundary
components of $R$ on its boundary.\\

(c) W.L.O.G. assume that $y$ is a boundary component of $R$, and
$z$ is a separating circle. Then, $z$ is a $(k, m+1)$ circle. By
using Lemma 2.11, it is easy to see that there exist $z' \in
\lambda([z])$ and an embedded pair of pants which has $x', z'$ and
a boundary component of $R$ on its boundary.\\

Hence, $\lambda$ ``sends'' an embedded pair of pants to an embedded pair of pants in case (i).\\

Case (ii): We can find a nonseparating circle $w$ and a
$(2,1)$-separating circle $t$ on $R$ such that $\{x, y, z, w\}$ is
pairwise disjoint and $x, y, z, w$ are on a genus 2 subsurface,
$R_1$, that $t$ bounds as shown in Figure 8. Let $P_1 = \{x, y, z,
w \}$. $P_1$ is a pair of pants decomposition of $R_1$. We can
complete $P_1 \cup \{t\}$  to a pants decomposition $P$ of $R$.\\

\begin{figure}[htb]
\label{picture4}
\begin{center}
\epsfig{file=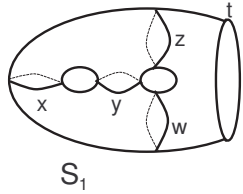,width=4.2cm} \caption{Nonseparating
circles, $x, y, z$, bounding a pair of pants}
\end{center}
\end{figure}

Let $R_2$ be the subsurface of $R$ of genus $g-2$ which is not
equal to $R_1$ and that comes from separation by $t$. By Lemma
\ref{division2}, there exist a $(2,1)$-separating circle $t' \in
\lambda([t])$ separating $R$ into two subsurfaces, $R_1', R_2'$,
of genus 2 and $g-2$ respectively such that
$\lambda(\mathcal{C}(R_1)) \subseteq \mathcal{C}(R_1')$ and
$\lambda(\mathcal{C}(R_2)) \subseteq \mathcal{C}(R_2')$. Since
$P_1 \subseteq R_1$, $\lambda([P_1]) \subseteq \mathcal{C}(R_1')$.
We can choose a set, $P_1'$, of pairwise disjoint representatives
of $\lambda([P_1])$ on $R_1'$. Then, $P_1' \cup \{t'\}$ is a pair
of pants decomposition of $R_1'$. We can choose a pairwise
disjoint representative set, $P'$, of $\lambda([P])$ containing
$P_1'$. $P'$ is a pair of pants decomposition of $R$. Let $x', y',
z', w' \in P_1'$ be the representatives of $\lambda([x]),
\lambda([y]), \lambda([z]), \lambda([w])$ respectively. Then,
since $t$ is adjacent to $z$ and $w$ w.r.t. $P$, $t'$ is adjacent
to $z'$ and $w'$ w.r.t. $P'$ by Lemma \ref{adjacent}. Then, since
$t' \cup z' \cup w' \subseteq R_1'$ and $t'$ is the boundary of
$R_1'$, there is an embedded pair of pants in $R_1'$ which has
$t', z', w'$ on its boundary. Since $z$ is a 4-curve in $P$, $z'$
is a 4-curve in $P'$. Since $z$ is adjacent to $x, y$ w.r.t. $P$,
$z'$ is adjacent to $x', y'$ w.r.t. $P'$. Since $z'$ is on the
boundary of a pair of pants which has $w', t'$ on its boundary,
and $z'$ is adjacent to $x', y'$, there is a pair of pants having
$x', y', z'$ on its boundary. So, $\lambda$ ``sends'' an embedded
pair of pants bounded by $x, y, z$ to an embedded pair of pants
bounded by $x', y', z'$ in this case too.\\

Case (iii): W.L.O.G assume that $z$ is a boundary component of $R$
and $x, y$ are nonseparating circles. Then, $([x], [y])$ is a
peripheral pair. Then, by Lemma 2.8, $(\lambda([x]),
\lambda([y]))$ is a peripheral pair. Let $x', y'$ be disjoint
representatives of $\lambda([x]), \lambda([y])$ respectively.
Then, there exists a pair of pants having $x', y'$ and a boundary
component of $R$ on its boundary. Hence, in this case also
$\lambda$ ``sends'' an embedded pair of pants to an embedded pair
of pants.\\

Now, assume that $P = (c_1, c_2, ..., c_{3g-3+p})$ is an ordered
pair of pants decomposition of $R$. Let $c_i' \in \lambda([c_i])$
such that the elements of $\{c_1', c_2', ..., c_{3g-3+p}'\}$ are
pairwise disjoint. Then, $P'=(c_1', c_2', ..., c_{3g-3+p}')$ is an
ordered pair of pants decomposition of $R$. Let $(B_1, B_2, ...,
B_{m})$ be an ordered set containing the connected components of
$R_P$. By the arguments given above, there is a corresponding,
``image'', ordered collection of pairs of pants $(B_1', B_2',...,
B_{m}')$. Nonembedded pairs of pants correspond to nonembedded
pairs of pants and embedded pairs of pants correspond to embedded
pairs of pants. Then, the proof of the lemma follows as in the
proof of Lemma 3.7 in [3].\end{proof}\\

\noindent {\bf Remark:} Let $\mathcal{E}$ be an ordered set of
vertices of $\mathcal{C}(R)$ having a pairwise disjoint
representative set $E$. Then, $E$ can be completed to an ordered
pair of pants decomposition, $P$, of $R$. We can choose an ordered
pairwise disjoint representative set, $P'$, of $\lambda([P])$ by
Lemma \ref{imageofpantsdecomp}. Let $E'$ be the elements of $P'$
which correspond to the elements of $E$. By Lemma \ref{top}, $P$
and $P'$ are topologically equivalent as ordered pants
decompositions. Hence, the set $E$ and $E'$ are topologically
equivalent. So, $\lambda$ gives a correspondence which preserves
topological equivalence on a set which has pairwise disjoint
representatives.\\


By using Lemma 2.10 and following the proof of Lemma 3.9 in [3],
we can prove the following lemma. This lemma will be used to see
some more properties of superinjective simplicial maps.\\

\begin{lemma}
\label{intone} Let $\lambda : \mathcal{C}(R) \rightarrow
\mathcal{C}(R)$ be a superinjective simplicial map. Let $\alpha$,
$\beta$ be two vertices of $\mathcal{C}(R)$. If $i(\alpha,
\beta)=1$, then $i(\lambda(\alpha), \lambda(\beta))=1$.
\end{lemma}

\section{Induced Map On Complex Of Arcs}

An arc $i$ on $R$ is called \textit{properly embedded} if
$\partial i \subseteq \partial R$ and $i$ is transversal to
$\partial R$. $i$ is called \textit{nontrivial} (or
\textit{essential}) if $i$ cannot be deformed into $\partial R$ in
such a way that the endpoints of $i$ stay in $\partial R$ during
the deformation. If $a$ and $b$ are two disjoint arcs connecting a
boundary component of $R$ to itself, $a$ and $b$ are called {\it
linked} if their end points alternate on the boundary component.
Otherwise, they are called {\it unlinked}.

\indent The \textit{complex of arcs}, $\mathcal{B}(R)$, on $R$ is
an abstract simplicial complex. Its vertices are the isotopy
classes of nontrivial properly embedded arcs $i$ in $R$. A set of
vertices forms a simplex if these vertices can be represented by
pairwise disjoint arcs.\\

\indent In this section, we assume that $\lambda : \mathcal{C}(R)
\rightarrow \mathcal{C}(R)$ is a superinjective simplicial map.
Let $\mathcal{V}(R)$ be the set of vertices of $\mathcal{B}(R)$.
We prove that $\lambda$ induces a map $\lambda_*: \mathcal{V}(R)
\rightarrow \mathcal{V}(R)$ with certain properties. Then we prove
that $\lambda_*$ extends to an injective simplicial map $\lambda_*
: \mathcal{B}(R) \rightarrow \mathcal{B}(R)$.\\

\indent The proofs of Lemma \ref{A} - \ref{D} are similar to the
proofs of the corresponding lemmas given in \cite{Ir}. So, we do
not prove these lemmas here. We only state them.

\begin{figure}[htb]
\begin{center} \epsfxsize=6.6in
\epsfbox{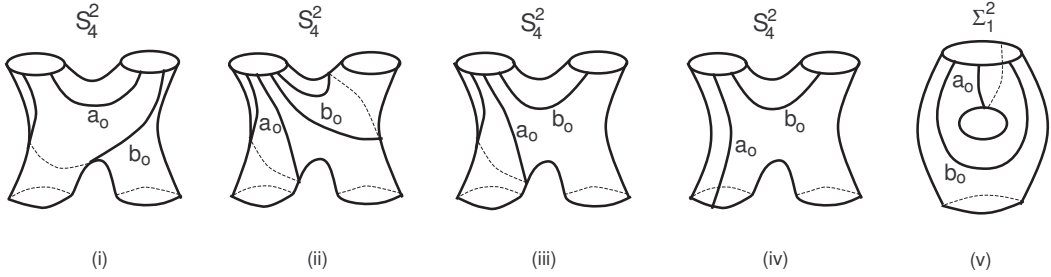} \caption{Disjoint arcs and neighborhoods}
\end{center}  \vspace{0.1cm}
\end{figure}

\begin{lemma}
\label{A} Let $a$ and $b$ be two disjoint arcs on $R$ connecting
two distinct boundary components, $\partial_i, \partial_j$, of
$R$. Let $N$ be a regular neighborhood of $a \cup b \cup
\partial_1 \cup \partial_2$ in $R$. Then, $(N, a, b) \cong (S_4
^2, a_o, b_o)$ where $S_4^2$ is a standard sphere with four holes
and $a_o, b_o$ are arcs as shown in Figure 9, (i).
\end{lemma}

\begin{lemma}
\label{B} Let $a$ and $b$ be two disjoint arcs which are unlinked,
connecting one boundary component $\partial_i$ of $R$ to itself.
Let $N$ be a regular neighborhood of $a \cup b \cup
\partial_i$ on $R$. Then, $(N,a, b) \cong (S_4^2, a_o, b_o)$
where $a_o, b_o$ are the arcs drawn on a standard sphere with four
holes, $S_4^2$, as shown in Figure 9, (ii).
\end{lemma}

\begin{lemma}
\label{C} Let $a$ and $b$ be two disjoint arcs on $R$ such that
$a$ connects one boundary component $\partial_i$ of $R$ to itself
for some $k=1,..., p$ and $b$ connects the boundary components
$\partial_i$ and $\partial_j$ of $R$, where $i \neq j$. Let $N$ be
a regular neighborhood of $a \cup b \cup \partial_i \cup
\partial_j$. Then, $(N, a, b) \cong (S_4^2, a_o, b_o)$ where $a_o,
b_o$ are the arcs drawn on a standard sphere with four holes,
$S_4^2$, as shown in Figure 9, (iii).
\end{lemma}

\begin{lemma}
\label{E} Let $a$ and $b$ be two disjoint arcs. Suppose that $a$
connects $\partial_i$ to $\partial_j$ and $b$ connects
$\partial_i$ to $\partial_k$ where $\partial_i$, $\partial_j$,
$\partial_k$ are three distinct boundary components. Let $N$ be a
regular neighborhood of $a \cup b \cup \partial_i \cup
\partial_j \cup \partial_k$. Then, $(N, a, b) \cong (S_4^2,
a_o, b_o)$ where $a_o, b_o$ are the arcs drawn on a standard
sphere with four holes, $S_4^2$, as shown in Figure 9,
(iv).\end{lemma} \vspace{0.1cm}

\begin{lemma}
\label{D} Let $a$ and $b$ be two disjoint, linked arcs connecting
one boundary component $\partial_i$ of $R$ to itself  for
$i=1,..., p$. Let $N$ be a regular neighborhood of $a \cup b \cup
\partial_i$. Then, $(N, a, b) \cong (\Sigma_1^2, a_o, b_o)$ where
$\Sigma_1^2$ is a standard surface of genus one with two boundary
components, and $a_o, b_o$ are as shown in Figure 9,
(v).\end{lemma} \vspace{0.1cm}

By using the following lemmas, we see some more properties of
$\lambda$.

\begin{lemma}
\label{horver} Let $\alpha$ and $\beta$ be two vertices in
$\mathcal{C}(R)$ which have representatives with geometric
intersection 2 and algebraic intersection 0 on $R$. Then,
$\lambda(\alpha)$ and $\lambda(\beta)$ have representatives with
geometric intersection 2 and algebraic intersection 0 on $R$.
\end{lemma}

\begin{figure}[htb]
\begin{center}
\epsfxsize=2.8in \epsfbox{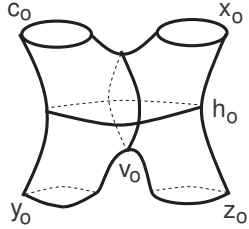} \caption{``Horizontal" and
``vertical" circles}
\end{center}
\end{figure}

\begin{proof} Let $h, v$ be representatives of $\alpha, \beta$
with geometric intersection 2 and algebraic intersection 0 on $R$.
W.L.O.G assume that $h$ and $v$ do not intersect the boundary
components of $R$. Let $N$ be a regular neighborhood of $h \cup v$
on $Int(R)$. Then, $N$ is a sphere with four boundary components.
Let $c, x, y, z$ be boundary components of $N$ such that there
exists a homeomorphism $\varphi : (N, c, x, y, z, h, v)$
$\rightarrow (N_o, c_o, x_o, y_o, z_o, h_o, v_o)$ where $N_o$ is a
standard sphere with four holes having $c_o, x_o, y_o, z_o$ on its
boundary and $h_o, v_o$ (horizontal, vertical) are two circles as
indicated in Figure 10. Since $h$ and $v$ have geometric
intersection 2 and algebraic intersection 0 on $R$, none of $c, x,
y, z$ bound a disk on $R$. If each of $c, x, y, z$ is an essential
circle on $R$, and $R$ is not a surface of genus two with two
boundary components, then the proof of the lemma follows from the
proof of Lemma 4.6 in [3], substituting Lemma 2.12 above for the
corresponding lemma in [3]. If each of $c, x, y, z$ is an
essential circle on $R$, and $R$ is a surface of genus two with
two boundary components, then the proof of the lemma follows from
the proof of Lemma 2.10. Assume that exactly one of $c, x, y, z$
is not essential. W.L.O.G assume that $c$ is not essential. Then,
since $N$ is a regular neighborhood in $Int(R)$, and $c$ does not
bound a disk on $R$, $c$ and a boundary component, say
$\partial_1$, of $R$ bound an annulus, $A$. Let $M= N \cup A$.
Then, $M$ is a regular neighborhood of $h \cup v$.\\

Let $A$= $\{x, y, z\}$. Any two elements in $A$ which are isotopic
in $R$ bound an annulus on $R$. Let $B$ be a set consisting of a
core from each annulus which is bounded by elements in $A$,
circles in $A$ which are not isotopic to any other circle in $A$,
and $v$. We can extend $B$ to a pants decomposition $P$ of $R$.
Since either $g = 2$ and $p \geq 2$ or $g \geq 3$ and $p > 0$,
there are at least four pairs of pants of $P$. Note that $\{v\}$
is a pair of pants decomposition of $M$. Each pair of pants of
this pants decomposition of $M$ is contained in exactly one pair
of pants in $P$. It is easy to see that there is a pair of pants
$Q$ of $P$ such that interior of $Q$ is disjoint from $M$, $Q$ has
at least one of $x, y, z$ as one of its boundary components and
all the boundary components of $Q$ are essential circles in $R$.
We will give the argument for the case where $Q$ has $y$ on its
boundary. The other cases follow by similar arguments.\\

Let $T$ be a regular neighborhood, in $Q$, of the boundary
components of $Q$ other than $y$. Let $t, w$ be the boundary
components of $T$ which are in the interior of $Q$. Then, $y, t,
w$ bound an embedded pair of pants $O$ in $Q$. Let $\tilde{M} = M
\cup O$. Then, we can extend $N_o$ to $\tilde{N_o}$ and find a
homeomorphism $\tilde{\varphi}: (\tilde{M}, \partial_1, x, y, z,
h, v, t, w) \rightarrow (\tilde{N_o}, c_o, x_o, y_o, z_o, h_o,
v_o, t_o, w_o)$, where $\tilde{N_o}$ is as shown in Figure 11.\\

Using Lemma \ref{top}, we can choose pairwise disjoint
representatives $x', y', z'$, $v', t', w'$ of $\lambda([x]),
\lambda([y]), \lambda([z]), \lambda([v]), \lambda([t]),
\lambda([w])$ respectively s.t. there exists a subsurface
$\tilde{M'}$ of $R$ and a homeomorphism $\chi : (\tilde{M'},
\partial_i, x', y', z', v', t', w')$ $ \rightarrow (\tilde{N_o},
c_o,$ $x_o, y_o, z_o, v_o, t_o, w_o)$ for some boundary component
$\partial_i$ of $R$. Since $i([h], [v]) \neq 0$ and $\lambda$ is
superinjective, we have, $i(\lambda([h]), \lambda([v])) \neq 0$.
Then, a representative $h'$ of $\lambda([h])$ can be chosen such
that $h'$ is transverse to $v'$, $h'$ doesn't intersect any of
$\partial_i, x', y', z'$, and $i(\lambda([h]), \lambda([v])) = |h'
\cap v'|$. Since $i(\lambda([h]), \lambda([v])) \neq 0$, $h'$
intersects $v'$. Hence, $h'$ is in the sphere with four holes
bounded by $\partial_i, x', y', z'$ in $\tilde{M'}$.\\

\begin{figure}
\begin{center}
\epsfxsize=2.8in \epsfbox{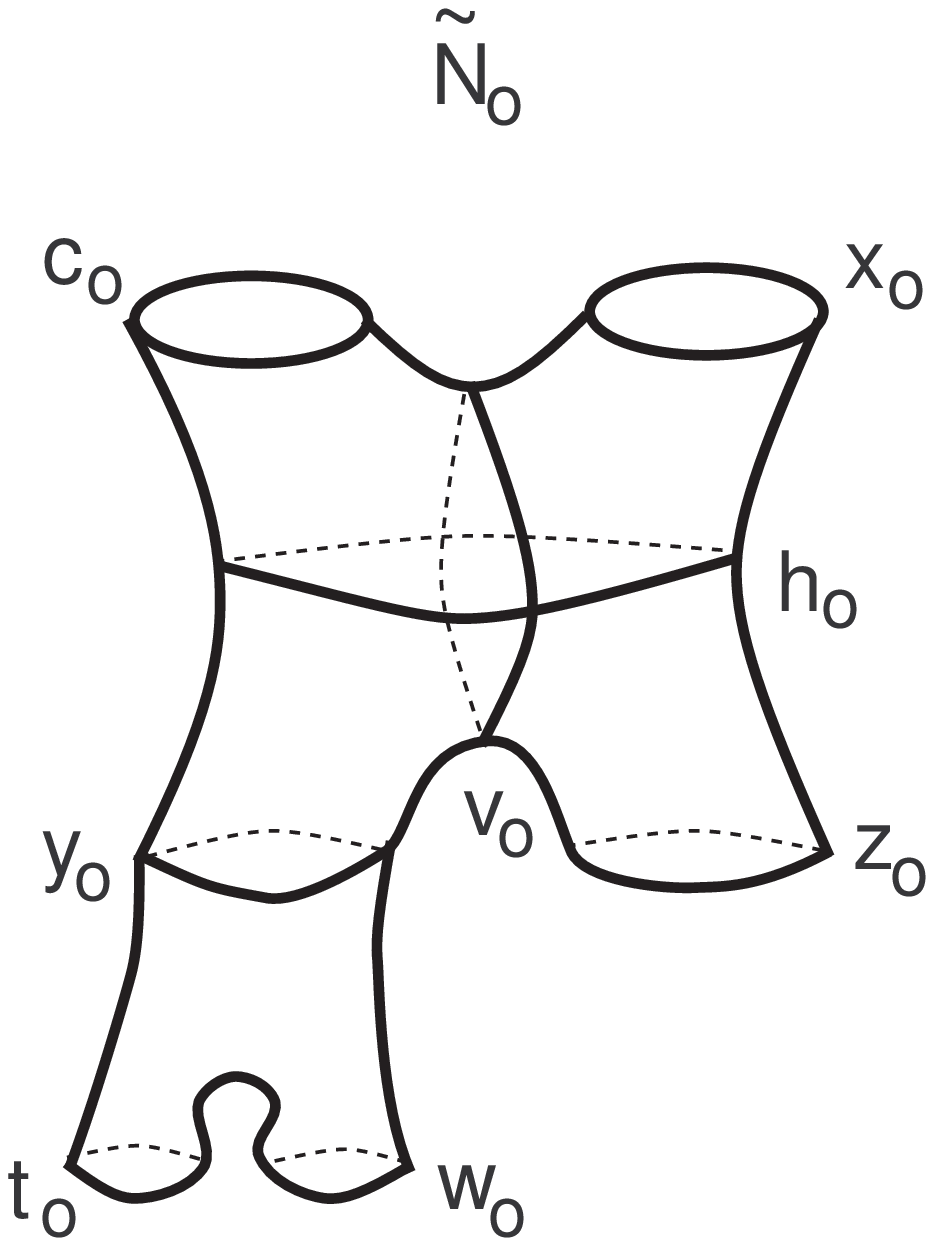} \caption{Sphere with five
holes}
\end{center}
\end{figure}

$\tilde{M}$ and $\tilde{M'}$ are spheres with five holes in $R$.
Since $x, z, t, w$ are essential circles in $R$, the essential
circles on $\tilde{M}$ are essential in $R$. Similarly, since $x',
z', t', w'$ are essential circles in $R$, the essential circles on
$\tilde{M'}$ are essential in $R$. Furthermore, we can identify
$\mathcal{C}(\tilde{M})$ and $\mathcal{C}(\tilde{M'})$ with two
subcomplexes of $\mathcal{C}(R)$ in such a way that the isotopy
class of an essential circle in $\tilde{M}$ or in $\tilde{M'}$ is
identified with the isotopy class of that circle in $R$. Now,
suppose that $\alpha$ is a vertex in $\mathcal{C}(\tilde{M})$.
Then, with this identification, $\alpha$ is a vertex in
$\mathcal{C}(R)$ and $\alpha$ has a representative in $\tilde{M}$.
Then, $i(\alpha, [x])=i(\alpha, [t])=i(\alpha, [w])= i(\alpha,
[z])=0$. Then there are two possibilities: (i) $\alpha = [v]$ or
$\alpha = [y]$, (ii) $i(\alpha, [v]) \neq 0$ or $i(\alpha, [y])
\neq 0$. Since $\lambda$ is injective, $\lambda(\alpha)$ is not
equal to any of $[x'], [t'], [w'], [z']$. Since $\lambda$ is
superinjective, $i(\lambda(\alpha), [x'])=i(\lambda(\alpha),
[t'])=i(\lambda(\alpha), [w'])= i(\lambda(\alpha), [z'])=0$. Then,
there are two possibilities: (i) $\lambda(\alpha) = [v']$ or
$\lambda(\alpha) = [y']$, (ii) $i(\lambda(\alpha), [v']) \neq 0$
or $i(\lambda(\alpha), [y']) \neq 0$. Then, a representative of
$\lambda(\alpha)$ can be chosen in $\tilde{M'}$. Hence, $\lambda$
maps the vertices of $\mathcal{C}(R)$ that have essential
representatives in $\tilde{M}$ to the vertices of $\mathcal{C}(R)$
that have essential representatives in $\tilde{M'}$, (i.e.
$\lambda$ maps $\mathcal{C}(\tilde{M}) \subseteq \mathcal{C}(R)$
to $\mathcal{C}(\tilde{M'}) \subseteq \mathcal{C}(R)$). Similarly,
$\lambda$ maps $\mathcal{C}(M) \subseteq \mathcal{C}(R)$ to
$\mathcal{C}(M') \subseteq \mathcal{C}(R)$.\\

It is easy to see that $\{[h], [v]\}$ is a simple pair in
$\tilde{M}$. Then, by Theorem \ref{korkmaz}, there exist vertices
$\gamma_1, \gamma_2, \gamma_3$ of $\mathcal{C}(\tilde{M})$ such
that $(\gamma_1, \gamma_2, [h], \gamma_3, [v])$ is a pentagon in
$\mathcal{C}(\tilde{M})$, $\gamma_1$ and $\gamma_3$ are
2-vertices, $\gamma_2$ is a 3-vertex, and $\{[h], \gamma_3\}$,
$\{[h], \gamma_2\}$, $\{[v], \gamma_3\}$ and $\{\gamma_1, \gamma_2
\}$ are codimension-zero simplices of $\mathcal{C}(\tilde{M})$.\\

Since $\lambda$ is superinjective and $x', t', w', z'$ are
essential circles, we can see that $(\lambda(\gamma_1),
\lambda(\gamma_2)$, $\lambda([h]), \lambda(\gamma_3)$,
$\lambda([v]))$ is a pentagon in $\mathcal{C}(\tilde{M'})$. By
Lemma \ref{top}, $\lambda(\gamma_1)$ and $\lambda(\gamma_3)$ are
2-vertices, and $\lambda(\gamma_2)$ is a 3-vertex in
$\mathcal{C}(\tilde{M'})$. Since $\lambda$ is an injective
simplicial map $\{\lambda([h]), \lambda(\gamma_3)\}$,
$\{\lambda([h]), \lambda(\gamma_2)\}$, $\{\lambda([v])$,
$\lambda(\gamma_3)\}$ and $\{\lambda(\gamma_1),
\lambda(\gamma_2)\}$ are codimension-zero simplices of
$\mathcal{C}(\tilde{M'})$. Then, by Theorem \ref{korkmaz},
$\{\lambda([h]), \lambda([v])\}$ is a simple pair in $\tilde{N'}$.
Since $\lambda([h])$ has a representative, $h'$, in $M'$, such
that $i(\lambda([h]), \lambda([v]) = |h' \cap v'|$ and
$\{\lambda([h]), \lambda([v])\}$ is a simple pair in $\tilde{M'}$,
there exists a homeomorphism $\chi : (M', \partial_i, x', y', z',
h', v') \rightarrow (N_o, c_o, x_o, y_o, z_o, h_o, v_o)$.\\

The proof of the lemma in the remaining cases, when $M$ has more
than one inessential boundary component is similar to the previous
case.\end{proof}

\begin{lemma}
\label{3} Let $c, x$ be curves which are either essential circles
on $R$ or some boundary components of $R$. Let $y, z, m, n$ be
essential circles on $R$ such that there exists a subsurface $N$
of $R$ and a homeomorphism $\varphi: (N, c, x, y, z, m, n)$ $
\rightarrow (N_o, c_o, x_o,$ $y_o, z_o, m_o, n_o)$ where $N_o$ is
a standard torus with two boundary components, $c_o, x_o$, and
$y_o, z_o$, $m_o, n_o$ are circles as shown in Figure 12 (i).
Then, there exist $c', x'$, two simple closed curves, and $y' \in
\lambda([y]), z' \in \lambda([z]), m' \in \lambda([m]), n' \in
\lambda([n]), N' \subseteq R$ and a homeomorphism $\chi: (N', c',
x', y', z'$, $m', n')$ $ \rightarrow (N_o, c_o,$ $x_o, y_o, z_o,
m_o, n_o)$.\end{lemma}

\begin{proof} If both $c$ and $x$ are essential circles,
then the proof follows from Lemma 4.7 in [3], substituting Lemma
2.12 above for the corresponding lemma in [3].\\

Since the genus of $R$ is at least 2, both of $c$ and $x$ cannot
be boundary components of $R$. W.L.O.G assume that $x$ is
essential and $c$ is a boundary component of $R$. We can complete
$\{x, y, z\}$ to a pair of pants decomposition, $P$, of $R$. Since
$\{y, z\}$ gives a pair of pants decomposition on $N$, by Lemma
\ref{top}, there exists a subsurface $N' \subseteq R$ which is
homeomorphic to $N_o$ and there are pairwise disjoint
representatives $x', y', z'$ of $\lambda([x]), \lambda([y]),
\lambda([z])$ respectively and a homeomorphism $\phi$ such that
$(N', \partial_i, x', y', z') \rightarrow _\phi (N_o, c_o, x_o,
y_o, z_o)$ for some $i \in \{1,...,k\}$. Then by Lemma
\ref{intone}, we have the following:\\

\begin{figure}[htb]
\begin{center}
\hspace{1.3cm} \epsfxsize=3.7in \epsfbox{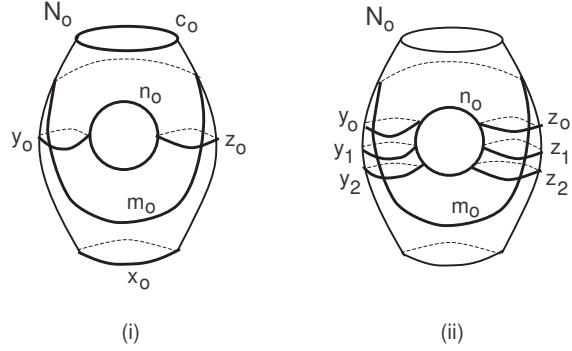}
\caption{Circles on torus with two boundary components}
\end{center}
\end{figure}

\noindent $i([m], [z])=1 \Rightarrow i(\lambda([m]),
\lambda([z]))=1$, $i([n], [z])=1 \Rightarrow i(\lambda([n]),
\lambda([z]))=1$, $i([m], [y])=1 \Rightarrow i(\lambda([m]),
\lambda([y]))=1$, $i([n], [y])=1 \Rightarrow i(\lambda([n]),
\lambda([y]))=1$,
$i([m], [x])=0 \Rightarrow i(\lambda([m]), \lambda([x]))=0$,
$i([n], [x])=0 \Rightarrow i(\lambda([n]), \lambda([x]))=0$,
$i([m], [n])=0 \Rightarrow i(\lambda([m]), \lambda([n]))=0$.\\

Let $c'=\partial_i$. There are representatives $m_1 \in
\lambda([m]), n' \in \lambda([n])$ such that $|m_1 \cap y'| = |m_1
\cap z'| = 1$, $|m_1 \cap c'| = |m_1 \cap x'| = 0$, $|n' \cap y'|
= |n' \cap z'| = 1, |n' \cap c'| = |n' \cap x'| = |m_1 \cap n'| =
0$ with all intersections transverse.\\

Since $\phi$ is a homeomorphism, we have $|\phi(m_1) \cap
\phi(n')|= 0, |\phi(n') \cap y_o| = 1$, $|\phi(n') \cap z_o| = 1,
|\phi(n') \cap c_o| = |\phi(n') \cap x_o| = 0 = |\phi(m_1) \cap
c_o| = |\phi(m_1) \cap x_o| = 0$, $|\phi(m_1) \cap y_o| =
|\phi(m_1) \cap z_o| =1$.\\

Let's choose parallel copies $y_1, y_2$ of $y_o$ and $z_1, z_2$ of
$z_o$ as shown in Figure 12 (ii) so that each of them has
transverse intersection one with $\phi(m_1)$ and $\phi(n')$. Let
$P_1, P_2$ be the pair of pants with boundary components $c_o,
y_o, z_o$, and $x_o, y_2, z_2$ respectively. Let $Q_1, Q_2, R_1,
R_2$ be the annulus with boundary components $\{y_o, y_1\}, \{y_1,
y_2\}$, $\{z_o, z_1\}$, $\{z_1, z_2\}$ respectively. By the
classification of isotopy classes of families of properly embedded
disjoint arcs in pairs of pants, $\phi(m_1) \cap P_1$, $\phi(m_1)
\cap P_2$, $\phi(n') \cap P_1$ and $\phi(n') \cap P_2$ can be
isotoped to the arcs $m_o \cap P_1, m_o \cap P_2, n_o \cap P_1,
n_o \cap P_2$ respectively. Let $\kappa : P_1 \times I \rightarrow
P_1$, $\tau : P_2 \times I \rightarrow P_2$ be such isotopies. By
a tapering argument, we can extend $\kappa $ and $\tau $ and get
$\tilde{\kappa }: (P_1 \cup Q_1 \cup R_1) \times I \rightarrow
(P_1 \cup Q_1 \cup R_1)$ and $\tilde{\tau }: (P_2 \cup Q_2 \cup
R_2) \times I \rightarrow (P_2 \cup Q_2 \cup R_2)$ so that
$\tilde{\kappa }_t$ is $id$ on $y_1 \cup z_1$ and $\tilde{\tau
}_t$ is $id$ on $y_1 \cup z_1$ for all $t \in I$. Then, by gluing
these extensions we get an isotopy $\vartheta$ on $N_o \times
I$.\\

By the classification of isotopy classes of arcs (relative to the
boundary) on an annulus, $\vartheta_1(\phi(n')) \cap (R_1 \cup
R_2)$ can be isotoped to $t_{z_o} ^k(n_o) \cap (R_1 \cup R_2)$ for
some $k \in \mathbb{Z}$. Let's call this isotopy $\mu$. Let
$\tilde{\mu }$ denote the extension by id to $N_o$. Similarly,
$\vartheta_1 (\phi(n')) \cap (Q_1 \cup Q_2)$ can be isotoped to
$t_{y_o} ^l (n_o) \cap (Q_1 \cup Q_2)$ for some $l \in
\mathbb{Z}$. Let's call this isotopy $\nu $. Let $\tilde{\nu }$
denote the extension by id to $N_o$. Then, ``gluing'' the two
isotopies $\tilde{\mu}$  and $\tilde{\nu}$, we get a new isotopy,
$\epsilon$, on $N_o$. Then we have, $t_{y_o} ^{-l} (t_{z_o}
^{-k}(\epsilon _1(\vartheta_1 (\phi(n'))))) = n_o$.
Clearly, $t_{y_o} ^{-l} \circ t_{z_o} ^{-k} \circ \epsilon_1 \circ
\vartheta _1$ fixes $c_o, x_o, y_o, z_o$. So, we get $t_{y_o}
^{-l} \circ t_{z_o}^{-k} \circ \epsilon_1 \circ \vartheta_1 \circ
\phi : (N', c', x', y', z', n') \rightarrow (N_o, c_o, x_o, y_o,
z_o, n_o)$. Let $\chi = t_{y_o} ^{-l} \circ t_{z_o}^{-k} \circ
\epsilon_1 \circ \vartheta_1 \circ \phi$. Then, because of the
intersection information we also have that $\chi(m_1)$ is isotopic
to either $m_o$ or $\hat{m}_o$ where $\hat{m}_o$ is the curve that
we get from $m_o$ by reflecting the picture in Figure 12 (i) about
the plane of the paper. Let $\rho$ be this reflection. We have
$\rho(m_o)=\hat{m}_o$. If $\chi(m_1)$ is isotopic to $m_o$, we let
$m'= \chi ^{-1}(m_o)$, and we get $\chi : (N', c', x', y', z', m',
n') \rightarrow (N_o, c_o, x_o, y_o, z_o, m_o, n_o)$. If
$\chi(m_1)$ is isotopic to $\hat{m}_o$, we let $m'= \chi
^{-1}(\hat{m}_o)$, and we get $\rho ^{-1} \circ \chi : (N', c',
x', y', z', m', n') \rightarrow (N_o, c_o, x_o, y_o, z_o, m_o,
n_o)$. This proves the lemma. \end{proof}\\

Let $i$ be an essential properly embedded arc on $R$. Let $A$ be a
boundary component of $R$ which has one end point of $i$ and $B$
be the boundary component of $R$ which has the other end point of
$i$. Let $N$ be a regular neighborhood of $i \cup A \cup B$ in
$R$. By Euler characteristic arguments, $N$ is a pair of pants.
The boundary components of $N$ are called {\it encoding circles of
$i$ on $R$}.
The set of isotopy classes of nontrivial encoding circles on $R$
is called the {\it encoding simplex}, $\Delta_i$, of $i$ (and of
$[i]$).\\

An essential properly embedded arc $i$ on $R$ is called {\it type
1} if it joins one boundary component $\partial_k$ of $R$ to
itself. It is called {\it type 1.1} if $\partial_k \cup i$ has a
regular neighborhood $N$ in $R$ which has only one circle on its
boundary which is inessential w.r.t. $R$. If $N$ has two circles
on its boundary which are inessential w.r.t. $R$, then $i$ is
called {\it type 1.2}. We call $i$ to be {\it type 2}, if it joins
two different boundary components of $R$ to each other. An element
$[i] \in \mathcal{V}(R)$ is called {\it type 1.1 (1.2, 2)} if it
has a type 1.1 (1.2, 2) representative. $i$ is called {\it
nonseparating} if its complement in $R$ is connected.\\

Let $\partial^1, \partial^2, ..., \partial^p$ be the boundary
components of $R$. We prove the following lemmas in order to show
that $\lambda$ induces a map $\lambda_* : \mathcal{V}(R)
\rightarrow \mathcal{V}(R)$ with certain properties.\\

\begin{lemma}
\label{boundary} Let $\partial_k \subseteq \partial R$ for some $k
\in \{1,..., p\}$. Then, there exists a unique $\partial^l \in
\partial R$ for some $l \in \{1,... ,p \}$ such that if $i$ is a
properly embedded essential arc on $R$ connecting $\partial_k$ to
itself, then there exists a properly embedded arc $j$ on $R$
connecting $\partial^l$ to itself such that $\lambda(\Delta_i) =
\Delta_j$.
\end{lemma}

\begin{proof} Assume that there are two boundary components,
$\partial^r$ and $\partial^t$ such that each of them satisfies the
hypothesis. Let $i$ be a properly embedded, essential,
nonseparating type 1 arc connecting $\partial_k$ to itself. Then,
there exist properly embedded arcs, $j_1$, connecting $\partial^r$
to itself, and $j_2$, connecting $\partial^t$ to itself, such that
$\lambda(\Delta_i) = \Delta_{j_1}$ and $\lambda(\Delta_i) =
\Delta_{j_2}$. Then, we have $\Delta_{j_1} = \Delta_{j_2}$. Note
that a properly embedded essential arc $i$ is type 1.1 iff
$\Delta_i$ has exactly 2 elements. Otherwise $\Delta_i$ has 1
element. Since $i$ is nonseparating type 1, it is type 1.1. Then,
since $\lambda(\Delta_i) = \Delta_{j_1}$ and $\lambda(\Delta_i) =
\Delta_{j_2}$ and $\lambda$ is injective, $j_1$ and $j_2$ are type
1.1. We can choose a pairwise disjoint representative set $\{a,
b\}$ of $\Delta_{j_1}$ on $R$. Since $\Delta_{j_1} =
\Delta_{j_2}$, $\{a, b\}$ is a pairwise disjoint representative
set for $\Delta_{j_2}$ on $R$. Then, $a, b$ and $\partial^r$ bound
a pair of pants, $P$, on $R$ containing an arc, $j_1'$, isotopic
to $j_1$. Similarly, $a, b,
\partial^t$ bound a pair of pants, $Q$, on $R$ containing an
arc, $j_2'$, isotopic to $j_2$. Let's cut $R$ along $a$ and $b$.
Then, $P$ is the connected component of $R_{a \cup b}$ containing
$\partial^r$ and $Q$ is the connected component of $R_{a \cup b}$
containing $\partial^t$. $P \neq Q$ since $\partial^t$ is not in
$P$ and $\partial^t$ is in $Q$. Then $P$ and $Q$ are distinct
connected components meeting along $a$ and $b$. Hence, $R$ is $P
\cup Q$, a torus with two holes which gives a contradiction since
the genus of $R$ is at least 2. So, only one boundary component of
$R$ can satisfy the hypothesis.\\

Since $i$ is nonseparating type 1, $\Delta_i$ contains two isotopy
classes of nontrivial circles in $R$. Let $P'$ be a pairwise
disjoint representative set of $\lambda([\Delta_i])$. Since the
genus of $R$ is at least 2, by the proof of Lemma \ref{top}, we
can see that $P'$ and a boundary component of $R$ bounds a unique
pair of pants $Q$, in $R$ which has only one inessential boundary
component. Let $\partial^{l(i)}$ be this inessential boundary
component. Let $j$ be an essential properly embedded arc
connecting $\partial^{l(i)}$ to itself in $Q$. Then, we have
$\lambda(\Delta_i) = \Delta_j$.\\

\begin{figure}[htb]
\begin{center}
\epsfxsize=2.65in \epsfbox{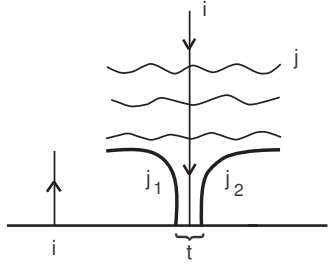} \caption{Splitting the arc
$j$ along the end of $i$}
\end{center}
\end{figure}

Now, to see that $\partial^{l(i)}$ is independent of the
nonseparating type 1 arc $i$ connecting $\partial_k$ to itself, we
prove the following claim:\\

Claim 1: If we start with two type 1 nonseparating arcs $i$ and
$j$ starting and ending on $\partial_k$, then $\partial^{l(i)} =
\partial^{l(j)}$.\\

Proof: Let $i$ and $j$ be nonseparating type 1 arcs connecting
$\partial_k$ to itself. W.L.O.G. we can assume that $i$ and $j$
have minimal intersection. First, we prove that there is a
sequence $j = r_0 \rightarrow r_1 \rightarrow ... \rightarrow
r_{n+1}=i$ of essential properly embedded nonseparating type 1
arcs joining $\partial_k$ to itself so that each consecutive pair
is disjoint, i.e. the isotopy classes of these arcs define a path
in $\mathcal{B} (R)$, between $i$ and $j$.\\

If $|i \cap j|=0$, then take $r_0=j$, $r_1=i$. Assume that $|i
\cap j|=m>0$. We orient $i$ and $j$ arbitrarily. Then, we define
two arcs in the following way: Start on the boundary component
$\partial_k$, on one side of the beginning point of $j$ and
continue along $j$ without intersecting $j$, till the last
intersection point along $i$. Then we would like to follow $i$,
without intersecting $j$, until we reach $\partial_k$. So, if we
are on the correct side of $j$ we do this; if not, we change our
starting side from the beginning and follow the construction. This
gives us an arc, say $j_1$. We define $j_2$, another arc, by
changing the orientation of $j$ and following the same
construction. It is easy to see that $j_1, j_2$ are disjoint
properly embedded arcs connecting $\partial_k$ to itself as $i$
and $j$ do. One can see that $j_1, j_2$ are essential arcs since
$i, j$ intersect minimally. In Figure 13, we show the beginning
and the end points of $i$, the essential intersections of $i, j$,
and $j_1, j_2$ near the end point of $i$ on $\partial_k$.\\

$|i \cap j_1| < m$, $|i \cap j_2| < m$ since we eliminated at
least one intersection with $i$. We also have $|j_1 \cap j| = |j_2
\cap j| =0$ since we never intersected $j$ in the construction.
Notice that $j_1$ and $j_2$ are not oriented, and $i$ is oriented.
It is easy to see from the construction that one of $j_1$ or $j_2$
has to be nonseparating type 1 arc, since $j$ is a nonseparating
type 1 arc.\\

Let $r_1 \in \{j_1, j_2\}$ and $r_1$ be nonseparating type 1. By
the construction, we get $|i \cap r_1| < m$, $|j \cap r_1|=0$.
Now, using $i$ and $r_1$ in place of $i$ and $j$ we can define a
new nonseparating type 1 arc $r_2$, with the properties $|i \cap
r_2| < |i \cap r_1|, |r_1 \cap r_2|=0$. By an inductive argument,
we get a sequence of nonseparating type 1 arcs such that every
consecutive pair is disjoint, $i=r_{n+1} \rightarrow r_n
\rightarrow r_{n-1} \rightarrow ... \rightarrow r_1 \rightarrow
r_0=j$. This gives us a special path in $\mathcal{B}(R)$ between
$i$ and $j$. By using Lemma \ref{B} and Lemma \ref{D}, we can see
a regular neighborhood of the union of each consecutive pair in
the sequence and the boundary component of $R$ that the arcs are
starting and ending at, and encoding circles of these consecutive
arcs. Then, by using the results of Lemma \ref{horver} and
\ref{3}, we can see that each pair of disjoint nonseparating type
1 arcs give us the same boundary component. Hence, by using the
sequence given above, we conclude that $i$ and $j$ give us the
same boundary component. This proves Claim 1.\\

\begin{figure}[htb]
\begin{center}
\epsfxsize=2.8in \epsfbox{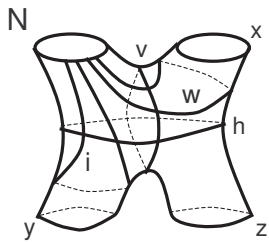} \caption{Unlinked arcs and
their encoding circles}
\end{center}
\end{figure}

Let $i_o$ be a properly embedded nonseparating type 1 arc on $R$
connecting $\partial_k$ to itself. Let $\partial^l =
\partial^{l(i_o)}$. If $i$ is a properly embedded, nonseparating
type 1 arc on $R$ connecting $\partial_k$ to itself, then by the
arguments given above we have $\partial^l = \partial^{l(i)}$, and
there exists a properly embedded arc $j$ on $R$ connecting
$\partial^l$ to itself such that $\lambda (\Delta_i) =
\Delta_j$.\\

Suppose that $i$ is a nontrivial properly embedded separating type
1 arc on $R$ connecting $\partial_k$ to itself. Then it is easy to
see that we can find an essential circle $v$ and a nonseparating
type 1 arc $w$, connecting $\partial_k$ to itself such that
$\partial_k \cup i \cup w$ has a regular neighborhood, $N$, which
is a sphere with four boundary components, as shown in Figure 14.
Let $x, y, z, h, v$ be as shown in the figure. Notice that $x$ and
$h$ are encoding circles of $w$. Since $w$ is a nonseparating type
1 arc, $x$ is essential. Then, by the proof of Lemma \ref{horver},
there exist essential simple closed curves $x' \in \lambda([x]),
h' \in \lambda([h]), v' \in \lambda([v]), N' \subseteq R$,
$\partial^{k'}, y', z'$, where $y', z'$ are boundary components of
$R$ if $y, z$ are inessential circles respectively and $y' \in
\lambda([y]), z' \in \lambda([z])$ is $y, z$ are essential circles
respectively and there exists a homeomorphism $\chi:(N',
\partial^{k'}, x', y', z', h',v') \rightarrow _\chi (N,
\partial_k, x, y, z, h, v)$.\\

Let $i'$ be an arc connecting $\partial^{k'}$ to itself in the
pair of pants determined by $\partial^{k'}, v', y'$ and $w'$ be an
arc connecting $\partial^{k'}$ to itself in the pair of pants
determined by $\partial^{k'}, x', h'$. We have $\lambda(\Delta_i)
= \Delta_{i'}$. Notice that since $w$ is a nonseparating type 1
arc, so is $w'$. Since $i'$ and $w'$ connect $\partial^{k'}$ to
itself, and $w'$ is a nonseparating type 1 arc, by using the
previous arguments, we see that $\partial^l = \partial^{k'}$. So,
the correspondence that we get on boundary components of $R$ using
nonseparating type 1 arcs is the same as the one that we get by
using separating type 1 arcs. Hence, $\partial^l$ is the boundary
component that we want.\end{proof}\\

We define a map $\sigma: \{ \partial_1, ..., \partial_p \}
\rightarrow \{ \partial^1, ..., \partial^p\}$ using the
correspondence which is given by Lemma \ref{boundary}.

\begin{lemma}
\label{vertex} Let $[i] \in \mathcal{V}(R)$. If $i$ connects
$\partial_k$ to $\partial_l$ on $R$ where $k, l \in \{1,..., p\}$,
then there exists a unique $[j] \in  \mathcal{V}(R)$ such that $j$
connects $\sigma(\partial_k)$ to $\sigma(\partial_l)$ and
$\lambda(\Delta_i) = \Delta_j$.
\end{lemma}

\begin{proof} Let $[i] \in \mathcal{V}(R)$ and let $i$ connect
$\partial_k$ to $\partial_l$ on $R$ where $k, l \in \{1,..., p\}$.
If $\partial_k = \partial_l$, there exists a nontrivial properly
embedded arc, $j$, connecting $\sigma(\partial_k)$ to itself such
that $\lambda(\Delta_i) = \Delta_j$ by Lemma 3.8. If $\partial_k
\neq \partial_l$, we can see the existence of $j$ as follows: Let
$a$ be a properly embedded nontrivial nonseparating arc which
connects $\partial_k$ to itself and let $b$ be a properly embedded
nontrivial nonseparating arc which connects $\partial_l$ to itself
such that $a, b, i$ are pairwise disjoint and they are on a
subsurface, $N$, which is a sphere with four boundary components,
as shown in Figure 15. Let $h, v, y, z$ be as shown in the figure.\\

\begin{figure}[htb]
\begin{center}
\epsfxsize=3.11in \epsfbox{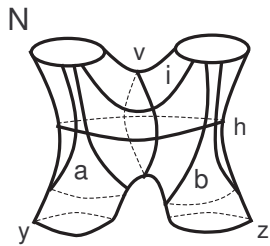} \caption{Arcs on sphere with
four holes}
\end{center}
\end{figure}

Since $a$ and $b$ are nonseparating, $y$ and $z$ are essential
circles. By Lemma \ref{horver}, there exists a subsurface $N'$,
representatives $h', v', y', z'$ in $\lambda([h]), \lambda([v]),
\lambda([y]),$ $ \lambda([z])$ respectively and two boundary
components, $\partial^r, \partial^t$ of $R$ and a homeomorphism
$\chi : (N, \partial_k, \partial_l, v, h, y, z) \rightarrow (N',
\partial^r, \partial^t, v', h', y', z')$. Then, by the proof of
Lemma \ref{boundary}, we see that $\partial^r =
\sigma(\partial_k),
\partial^t = \sigma(\partial_l)$. Let $j$ be a properly embedded
arc connecting $\sigma(\partial_k)$ to $\sigma(\partial_l)$ in the
pair of pants bounded by $\sigma(\partial_k), \sigma(\partial_l)$
and $h'$ . Then, we have $\lambda(\Delta_i) = \Delta_j$.\\

Now, let $e$ be an essential properly embedded arc in $R$ such
that $e$ connects $\sigma(\partial_k)$ to $\sigma(\partial_l)$ and
$\lambda(\Delta_i) = \Delta_e$. Then, we have $\Delta_e = \Delta_j
= \lambda(\Delta_i)$. Let $Q$ be a regular neighborhood of $e \cup
\sigma(\partial_k) \cup \sigma(\partial_l)$. Since $\Delta_e =
\Delta_j$, there is a properly embedded arc $j_1$ isotopic to $j$
in $Q$. Then, since both $j_1$ and $e$ connect the same boundary
components in this pair of pants, they are isotopic. Then, $[j] =
[e]$. Hence, $[j]$ is the unique isotopy class in $R$ such that
$j$ connects $\sigma(\partial_k)$ to $\sigma(\partial_l)$ and
$\lambda(\Delta_i) = \Delta_j$. \end{proof}\\

$\lambda$ induces a unique map $\lambda_* : \mathcal{V}(R)
\rightarrow \mathcal{V}(R)$ such that if $[i] \in \mathcal{V}(R)$
then $\lambda_*([i])$ is the unique isotopy class corresponding to
$[i]$ where the correspondence is given by Lemma \ref{vertex}.
Using the results of the following lemmas, we will prove that
$\lambda_*$ extends to an injective simplicial map on
$\mathcal{B}(R)$.\\

\begin{lemma} $\lambda_* :\mathcal{V}(R) \rightarrow \mathcal{V}(R)$
extends to a simplicial map $\lambda_* : \mathcal{B}(R)
\rightarrow \mathcal{B}(R)$.
\end{lemma}

\begin{proof} It is enough to prove that if two distinct isotopy
classes of essential properly embedded arcs on $R$ have disjoint
representatives, then their images under $\lambda_*$ have disjoint
representatives. Let $a, b$ be two disjoint representatives of two
distinct classes in $\mathcal{V}(R)$. Let $\partial_1, ...,
\partial_p$ be the boundary components of $R$. We consider
the following cases:\\

\begin{figure}[htb]
\begin{center}
\epsfxsize=5.1in \epsfbox{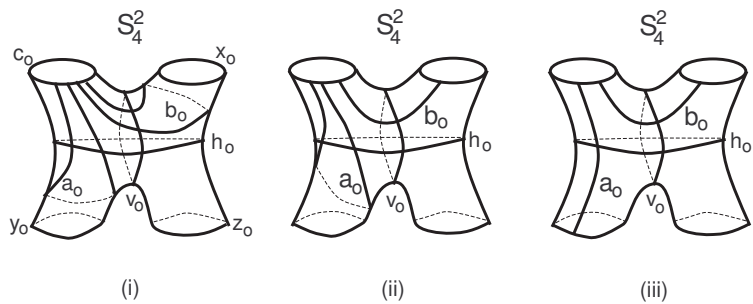} \caption{Arcs and their
encoding circles}
\end{center}
\end{figure}

Case 1: Assume that there is no common boundary among the
boundaries that $a$ and $b$ connect. Let $\partial_i, \partial_j$
be the boundary components that $a$ starts and ends, and let
$\partial_m, \partial_n$ be the boundary components that $b$
starts and ends where $i,j= 1, ..., p$ and $m,n= 1,..., p$ and
$\{m,n\} \cap \{i,j\} = \emptyset$. Then, since $a \cup
\partial_i \cup \partial_j$ is disjoint from $b \cup \partial_m
\cup \partial_n$, we can find disjoint regular neighborhoods,
$N_1$, of $a \cup \partial_i \cup \partial_j$ and $N_2$ of $b \cup
\partial_m \cup \partial_n$ on $R$, which give us two disjoint
pair of pants. Then, by using Lemma \ref{top}, and the definition
of $\lambda_*$, it is easy to see that the corresponding arcs
(images) will have disjoint representatives.\\

Case 2: Assume that $a, b$ are unlinked, connecting $\partial_i$
to itself for some $i=1,..., p$. By Lemma \ref{B}, there is a
homeomorphism $\phi$ such that $(S^2_4, a_o, b_o) \cong_\phi (N,
a, b)$ where $N$ is a regular neighborhood of $a \cup b \cup
\partial_i$ in $R$ and $a_o, b_o$ are as shown in Figure 16 (i).
Then, by using Lemma \ref{horver} and the definition of
$\lambda_*$, we see that images have disjoint representatives.\\

Case 3: Assume that $a$ connects one boundary component
$\partial_i$ to itself for some $i=1,..., p$, and $b$ connects
$\partial_i$ to $\partial_k$ for some $k \neq i$. Then, by Lemma
\ref{C}, there is a homeomorphism $\phi$ such that $(S^2_4, a_o,
b_o) \cong_\phi (N, a, b)$ where $N$ is a regular neighborhood of
$a \cup b \cup \partial_i \cup \partial_k$ in $R$ and $a_o, b_o$
are as shown in Figure 16 (ii). Then, by using Lemma \ref{horver}
and the definition of $\lambda_*$, we see that the images have
disjoint representatives.\\

Case 4: Assume that $a$ connects $\partial_i$ to $\partial_j$ and
$b$ connects $\partial_i$ to $\partial_k$ where $k \neq i, i \neq
j, j \neq k$. Then, by Lemma \ref{E}, there is a homeomorphism
$\phi$ such that $(S^2_4, a_o, b_o) \cong_\phi (N, a, b)$ where
$N$ is a regular neighborhood of $a \cup b \cup \partial_i \cup
\partial_k \cup \partial_j$ in $R$ and $a_o, b_o$ are as shown in
Figure 16 (iii). Then, by using Lemma \ref{horver} and the
definition of $\lambda_*$, as in the previous cases we see that
the images have disjoint representatives.\\

\begin{figure}[htb]
\begin{center}
\epsfxsize=4.2in \epsfbox{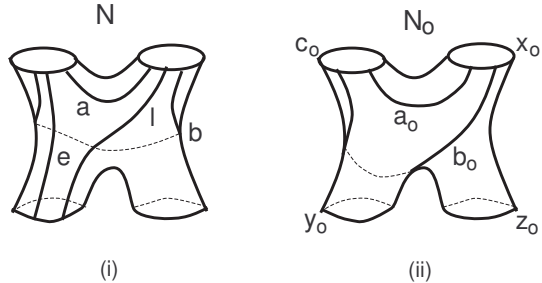} \caption{Arcs on sphere with
four holes}
\end{center}
\end{figure}

Case 5: Assume that each of $a$ and $b$ connect two distinct
boundary components, say $\partial_i, \partial_j$, of $R$.\\

We will first prove this case when $R$ has at least 3 boundary
components. Let $\partial_k$ be a boundary component different
from $\partial_i$ and $\partial_j$. Let $e$ and $l$ be disjoint
properly embedded arcs which are disjoint from $a, b$ such that
$e$ connects $\partial_i$ to $\partial_k$, $l$ connects
$\partial_j$ to $\partial_k$ and $a, e, l$ are in a subsurface
$N$, which is a sphere with 4 holes, of $R$ as shown in Figure 17
(i). Then, by applying the result of case 4 to each pair in $\{a,
e, l\}$, we can find disjoint representatives $a_1, e_1$ of
$\lambda_*([a]), \lambda_*([e])$ respectively, disjoint
representatives $a_2, l_1$ of $\lambda_*([a]), \lambda_*([l])$
respectively, disjoint representatives, $e_2, l_2$ of
$\lambda_*([e]), \lambda_*([l])$ respectively. Then, by using
Lemma \ref{E}, and Lemma \ref{horver} we can choose disjoint
representatives $a', e', l'$ of $\lambda_*([a]), \lambda_*([e]),
\lambda_*([l])$ respectively, a subsurface $N' \subseteq R$ and a
homeomorphism $(N', a', e', l') \rightarrow (N, a, e, l)$ where
$N, a, e, l$ are as shown in Figure 17 (i). Since $b$ and $e$ are
disjoint and $b$ connects $\partial_i$ to $\partial_j$ and $e$
connects $\partial_i$ to $\partial_k$, by using case 4, we can
choose a representative $b_1$ of $\lambda_*([b])$ which is
disjoint from $e'$. Similarly, we can choose a representative
$b_2$ of $\lambda_*([b])$ which is disjoint from $l'$. Then, since
$e'$ and $l'$ are disjoint, we can choose a representative $b_3$
of $\lambda_*([b])$ which is disjoint from $e' \cup l'$. Then, it
is easy to see that $\lambda([b])$ has a representative which is
disjoint from $a'$.\\

Now, assume that $R$ has exactly two boundary components,
$\partial_i, \partial_j$. By Lemma \ref{A}, there is a
homeomorphism $\phi: (N_o, a_o, b_o, c_o, x_o, y_o, z_o)
\rightarrow (N, a, b,$ $\partial_i, \partial_j, y, z)$ where $N$
is a regular neighborhood of $a \cup b \cup
\partial_i \cup \partial_j$ in $R$ and $a_o, b_o, c_o, x_o,
y_o, z_o$ are as shown in Figure 17 (ii). Since $R$ has exactly
two boundary components, $y$ and $z$ are essential circles in $R$.
Then, using Lemma \ref{top}, we can choose pairwise disjoint
representatives $y', z'$ of $\lambda([y]), \lambda([z])$
respectively, boundary components $\partial^k, \partial^l$ and a
subsurface $N' \subseteq R$, and a homeomorphism $\chi : (N',
\partial^k, \partial^l, y', z') \rightarrow (S_4^2, c_o, x_o, y_o,
z_o)$.\\

$N$ and $N'$ are spheres with four holes in $R$. Since $y, z$ are
essential circles in $R$, the essential circles on $N$ are
essential in $R$. Similarly, since $y', z'$ are essential circles
in $R$, the essential circles on $N'$ are essential in $R$.
Furthermore, we can identify $\mathcal{C}(N)$ and
$\mathcal{C}(N')$ with two subcomplexes of $\mathcal{C}(R)$ in
such a way that the isotopy class of an essential circle in $N$ or
in $N'$ is identified with the isotopy class of that circle in
$R$. Now, suppose that $\alpha$ is a vertex in $\mathcal{C}(N)$.
Then, with this identification, $\alpha$ is a vertex in
$\mathcal{C}(R)$ and $\alpha$ has a representative in $N$. Then as
in the proof of Lemma \ref{horver}, we can see that $\lambda$ maps
$\mathcal{C}(N) \subseteq \mathcal{C}(R)$ to $\mathcal{C}(N')
\subseteq \mathcal{C}(R)$. Then, since $N$ and $N'$ has four
boundary components, we can apply the arguments given in the proof
of the first part to see that there are disjoint representatives
of $\lambda_*([a])$ and $\lambda_*([b])$ in $N'$. This proves case
5.\\

\begin{figure}[htb]
\begin{center}
\epsfxsize=2.75in \epsfbox{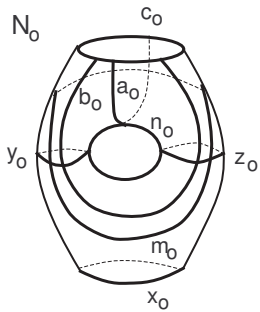} \caption{Linked arcs and
their encoding circles} \label{Figure18}
\end{center}
\end{figure}

Case 6: Assume that $a, b$ are linked, connecting $\partial_i$ to
itself for some $i=1, ...,p$. By Lemma \ref{D}, there is a
homeomorphism $\phi:(N_o, a_o, b_o) \rightarrow (N, a, b)$ where
$N$ is a regular neighborhood of $\partial_i \cup a \cup b$ in $R$
and $N_o, a_o, b_o$ are as in Figure \ref{Figure18}. Since $y_o,
z_o, c_o$ and $m_o, n_o, c_o$ are the boundary components of
regular neighborhoods of $a_o \cup c_o$ and $b_o \cup c_o$ on
$N_o$ respectively, $\phi(y_o), \phi(z_o)$ and $\phi(m_o),
\phi(n_o)$ are encoding circles for $a$ and $b$ on $R$
respectively. We have $(N, \phi(c_o), \phi(x_o), \phi(y_o),
\phi(z_o), \phi(m_o), \phi(n_o))$ $ \cong$ $(N_o, c_o,$ $x_o, y_o,
z_o, m_o, n_o)$. Since the genus of $R$ is at least 2, $\phi(x_o)$
is an essential circle on $R$. By Lemma \ref{3}, there exists
$\partial^j$, $x_o' \in \lambda (\phi(x_o))$, $y_o' \in
\lambda(\phi(y_o)), z_o' \in \lambda(\phi(z_o)), m_o' \in \lambda
(\phi(m_o)), n_o' \in \lambda (\phi(n_o)$, $N' \subseteq R$ and a
homeomorphism $\chi : (N_o, c_o, x_o, y_o, z_o, m_o$, $n_o)
\rightarrow (N', \partial^j, x_o', y_o', z_o', m_o', n_o')$. Since
$\phi(y_o), \phi(z_o)$ and $\phi(m_o), \phi(n_o)$ are encoding
circles for $a$ and $b$ on $R$ respectively, $y_o', z_o'$ and
$m_o', n_o'$ are encoding circles for $\lambda_*([a]),
\lambda_*([b])$ on $R$ respectively. Existence of $\chi$ shows
that $\lambda_*([a]), \lambda_*([b])$ have disjoint
representatives. $\chi(a_o)$ and $\chi(b_o)$ are disjoint
representatives for $\lambda_*([a]), \lambda_*([b])$
respectively.\\

We have shown in all the cases that if two vertices have disjoint
representatives, then $\lambda_*$ sends them to two vertices which
have disjoint representatives. Hence, $\lambda_*$ extends to a
simplicial map $\lambda_*: \mathcal{B}(R) \rightarrow
\mathcal{B}(R)$.\end{proof}

\begin{lemma}
\label{inj2} Let $\lambda : \mathcal{C}(R) \rightarrow
\mathcal{C}(R)$ be a superinjective simplicial map. Then,
$\lambda_* : \mathcal{B}(R) \rightarrow \mathcal{B}(R)$ is
injective.\end{lemma}

\begin{proof} It is enough to prove that $\lambda_*$ is injective
on the vertex set, $\mathcal{V}(R)$. Let $[i], [j] \in
\mathcal{V}(R)$ such that $\lambda_*([i])= \lambda_*([j])=[k]$.
Then, by the definition of $\lambda_*$, the type of [i] and [j]
are the same. Assume they are both type 1.1. Let $\{[x], [y]\}$
and $\{[z], [t]\}$ be the encoding simplices for $[i]$ and $[j]$
respectively. Then, $\{\lambda([x]), \lambda([y])\}$ and
$\{\lambda([z]), \lambda([t])\}$ are encoding simplices of $[k]$.
So, \{$\lambda ([x]), \lambda ([y])\} = \{\lambda ([z]), \lambda
([t])\}$. Then, since $\lambda$ is injective, we get $\{[x], [y]\}
= \{[z], [t]\}$. This implies $[i]=[j]$. The other cases can be
proven similarly to the first case by using the injectivity of
$\lambda$.\end{proof}\\

The following lemma can be proven similar to the proof of Lemma
4.13 in \cite{Ir}, which uses the Connectivity Theorem for
Elementary Moves of Mosher, \cite{Mos}, appropriately restated for
surfaces with boundaries. We will only state this lemma here.\\

\begin{lemma}
\label{fixing} If an injection $\mu : \mathcal{B}(R) \rightarrow
\mathcal{B}(R)$ agrees with $h_* : \mathcal{B}(R) \rightarrow
\mathcal{B}(R)$ on a top dimensional simplex, where $h_*$ is
induced by a homeomorphism $h: R \rightarrow R$, then $\mu$ agrees
with $h_*$ on $\mathcal{B}(R)$.\end{lemma}

\noindent {\bf Notation:} A homeomorphism $g: R \rightarrow R$
induces a map $g_\#: \mathcal{C}(R) \rightarrow \mathcal{C}(R)$,
where $g_\# =[g]$ and $g_\#$ induces a map $g_*: \mathcal{B}(R)
\rightarrow \mathcal{B}(R)$ in a similar way as $\lambda$ induces
$\lambda_*$.\\

\noindent {\bf Remark:} We have proven that $\lambda$ is an
injective simplicial map which preserves the geometric
intersection 0 and 1. If the number of boundary components of $R$
is at least 2, using these properties and following N.V.Ivanov's
proof of his Theorem 1.1 \cite{Iv1}, it can be seen that
$\lambda_*$ agrees with a map, $h_\#$, induced by a homeomorphism
$h: R \rightarrow R$ on a top dimensional simplex in
$\mathcal{B}(R)$. Then, by Lemma \ref{fixing}, it agrees with
$h_\#$ on $\mathcal{B}(R)$. Then, it is easy to see that $\lambda$
agrees with a map, $h_*$ on $\mathcal{C}(R)$.\\

If the number of boundary components of $R$ is exactly 1, we cut
$R$ along a nonseparating simple closed curve, $c$ on $R$. Let
$R_c$ be this cut surface and let $c_+, c_-$ be the two boundary
components of $R_c$ which comes from cutting $R$ along $c$.  Then,
considering how $\lambda$ induced $\lambda_*$ and using the
techniques of \cite{Ir}, it is easy to see that $\lambda$ induces
a superinjective simplicial map $\lambda_c: \mathcal{C}(R_c)
\rightarrow \mathcal{C}(R_d)$ and an injective simplicial map
$(\lambda_{c})_*: \mathcal{B}(R_c) \rightarrow \mathcal{B}(R_d)$,
where $\lambda([c]) = [d]$. Then, since $R_c$ and $R_d$ have 3
boundary components, by adapting the arguments given in the
paragraph above and using Lemma \ref{top}, we see that $\lambda_c$
agrees with a map $(g_c)_\#$ on $\mathcal{C}(R_c)$, where
$(g_c)_\#$ is induced by a homeomorphism $g_c: R_c \rightarrow
R_d$ such that $g_c(\{c_+, c_-\})= \{d_+, d_-\}$. We can do this
argument for any nonseparating simple closed curve on $R$. Then,
to see that $\lambda$ agrees with a map, $h_\#$, which is induced
by a homeomorphism $h: R \rightarrow R$, on $\mathcal{C}(R)$ we
use the following lemma.

\begin{lemma}
\label{imp} Let $\lambda: \mathcal{C}(R) \rightarrow
\mathcal{C}(R)$ be a superinjective simplicial map. Assume that
for any nonseparating simple closed curve $c$ on $R$, $\lambda_c$
agrees with a map, $(g_c)_\# : \mathcal{C}(R_c) \rightarrow
\mathcal{C}(R_d)$, which is induced by a homeomorphism $g_c: R_c
\rightarrow R_d$ where $g_c(\{c_+, c_-\})= \{d_+, d_-\}$ and
$\lambda([c]) = [d]$. Then, $\lambda$ agrees with a map $h_\#:
\mathcal{C}(R) \rightarrow \mathcal{C}(R)$ which is induced by a
homeomorphism $h: R \rightarrow R$.
\end{lemma}

\begin{proof} Let $c$ be a nonseparating simple closed curve
and $(g_c)_\# : \mathcal{C}(R_c) \rightarrow \mathcal{C}(R_d)$ be
a simplicial map induced by a homeomorphism $g_c: R_c \rightarrow
R_d$ where $g_c(\{c_+, c_-\})= \{d_+, d_-\}$ and $\lambda([c]) =
[d]$ such that $\lambda_c$ agrees with $(g_c)_\#$ on
$\mathcal{C}(R_c)$. Let $g$ be a homeomorphism of $R$ which cuts
to a homeomorphism $R_c \rightarrow R_d$ which is isotopic to
$g_c$. Then each homeomorphism of $R$ which cut to a homeomorphism
$R_c \rightarrow R_d$ which is isotopic to $g_c$, is isotopic to
an element in the set $\{gt_c^n, n \in \mathbb{Z}\}$, [1]. It is
easy to see that $\lambda_c$ agrees with $((gt_c^n)_c)_\#$ on
$\mathcal{C}(R_c)$ for all $n \in \mathbb{Z}$.\\

Let $w$ be a simple closed curve which is dual to $c$ (i.e. $w$
intersects $c$ transversely once and there is no other
intersection). Let $P$ be a regular neighborhood of $c \cup w$.
Then, $P$ is a genus one surface with one boundary component. Let
$y$ be the boundary component of $P$. We have $i([c], [y]) = 0$,
$i([w], [y]) = 0$, and $i([c], [w]) = 1$. Then, since $\lambda$ is
superinjective we have $i(\lambda([c]), \lambda([y])) = 0$,
$i(\lambda([w]), \lambda([y])) = 0$ and $i(\lambda([c]),
\lambda([w])) = 1$.\\

Let $Q$ be the genus one subsurface with one boundary component of
$R$ which has $g(y)$ as its boundary. Then, it is easy to see that
$g(c)=d \subseteq Q$, $g(w) \subseteq Q$ and $g(w)$ is dual to
$d$, since $w$ is dual to $c$.\\

Since $[y] \in \mathcal{C}(R_c)$, $\lambda([y])= g_\#([y]) =
[g(y)]$. We also have $[d]= \lambda([c])= g_\#([c])$. Since
$i(\lambda([c]), \lambda([y])) = 0$, $i(\lambda([w]),
\lambda([y])) = 0$ and $i(\lambda([c]), \lambda([w])) = 1$, and $d
\in \lambda([c])$, $g(y) \in \lambda([y])$ and $d$ and $g(y)$ are
disjoint, we can choose a simple closed curve $w' \in
\lambda([w])$ such that $w'$ is in $Int(Q)$ and dual to $d$. Then,
$g^{-1}(w')$ is dual to $c$ and $g^{-1}(w')$ is in $Int(P)$. Then,
since both of $w$ and $g^{-1}(w')$ are dual to $c$ in $Int(P)$,
there exists $m_c \in \mathbb{Z}$ such that $t_c^{m_c}([w]) =
[g^{-1}(w')]$. Then, $gt_c^{m_c} ([w]) = [w']$. Since
$\lambda([w]) = [w']$, $gt_c^{m_c}$ agrees with $\lambda$ on
$[w]$. We can identify $\mathcal{C}(R_c)$ with a subcomplex,
$L_c$, of $\mathcal{C}(R)$. Let $D_c$ be the set of isotopy
classes of simple closed curves which are dual to $c$ on $R$.\\

Claim 1: $(gt_c^{m_c})_\#$ agrees with $\lambda$ on $\{[c]\} \cup
L_c \cup D_c$.\\

Proof: It is clear that $(gt_c^{m_c})_\# ([c])= \lambda ([c]) =
[d]$. Since $(g_c)_\#$ agrees with $\lambda_c$ on
$\mathcal{C}(R_c)$, $(gt_c^{m_c})_\#$ agrees with $\lambda$ on
$L_c$.\\

We have seen that $gt_c^{m_c}$ agrees with $\lambda$ on $[w]$. Let
$w_1$ be a simple closed curve which is disjoint from $w$ and dual
to $c$. As we described before, there exists $\tilde{m_c} \in
\mathbb{Z}$ such that $\lambda$ agrees with $gt_c^{\tilde{m_c}}$
on $[w_1]$. Since $w$ and $w_1$ are disjoint, $i(\lambda([w]),
\lambda([w_1])) = 0$. If $m_c \neq \tilde{m_c}$ then
$i((gt_c^{m_c})(w) , gt_c^{\tilde{m_c}}(w_1)) \neq 0$ (since both
$w$ and $w_1$ are dual to $c$). Then, since $\lambda([w])=
(gt_c^{m_c})([w])$ and $\lambda([w_1])=
(gt_c^{\tilde{m_c}})([w_1])$, we would get $i(\lambda([w]),
\lambda([w_1])) \neq 0$, which gives a contradiction. Therefore,
$m_c = \tilde{m_c}$. Then, we see that $(gt_c^{m_c})_\#$ agrees
with $\lambda$ on $\{[c]\} \cup L_c \cup \{ [w] \cup [w_1] \}$.\\

Given any curve $t$ which is dual to $c$, by using similar
techniques as in Lemma 3.8, we can find a sequence of dual curves
to $c$, connecting $d$ to $t$, such that each consecutive pair is
disjoint, i.e. the isotopy classes of these curves define a path
in $\mathcal{C} (R)$, between $d$ and $t$. Then using the argument
above and the sequence, we conclude that $(gt_c^{m_c})_\#$ agrees
with $\lambda$ on $D_c$. Hence, $(gt_c^{m_c})_\#$ agrees with
$\lambda$ on $\{[c]\} \cup L_c \cup D_c$. This proves claim 1. Let
$h_c = gt_c^{m_c}$. We have that $(h_c)_\#$ agrees with $\lambda$
on $\{[c]\} \cup L_c \cup D_c$.\\

Claim 2: Let $v$ be a nonseparating simple closed curve on $R$.
Then, $(h_c)_\# = (h_v)_\# = \lambda$ on $\mathcal{C}(R)$.\\

Proof: Since the genus of $R$ is at least 2, we can find a
sequence of nonseparating simple closed curves connecting $c$ to
$v$ such that each consecutive pair is disjoint. Refining this
sequence, we can get a sequence $c \rightarrow c_1 \rightarrow ...
\rightarrow c_{n}=v$ of nonseparating simple closed curves
connecting $c$ to $v$ such that each consecutive pair in this
sequence is simultaneously nonseparating.\\

\begin{figure}
\begin{center}
\epsfxsize=2.3in \epsfbox{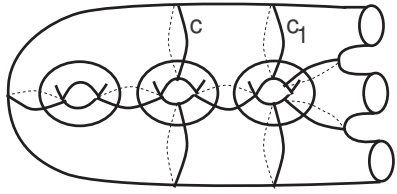} \caption{A configuration of
circles}
\end{center}
\end{figure}

Let's consider the first consecutive pair in the sequence, $c,
c_1$. Since $\{c, c_1\}$ is simultaneously nonseparating, it can
be completed to a set $G$, (shown in Figure 19, for $g=3, p=3$),
see \cite{IMc}, such that the isotopy classes of Dehn twists about
the elements of this set generate $PMod_R$ and all the curves in
this set are (i) either disjoint from $c$ or dual to $c$, and (ii)
either disjoint from $c_1$ or dual to $c_1$. Then, since all the
curves in $G$ are either disjoint from $c$ or dual to $c$, by
Claim 1 we have that $(h_c)_\# ([x]) = \lambda ([x])$ for every $x
\in G$. Similarly, since all the curves in $G$ are either disjoint
from $c_1$ or dual to $c_1$, by Claim 1 we have $(h_{c_1})_\#
([x]) = \lambda ([x])$ for every $x \in G$. Hence, $(h_c)_\#
([x])= \lambda ([x]) = (h_{c_1})_\# ([x])$ for every $x \in G$.
Then, $(h_{c})_\# = (h_{c_1})_\#$ since $(h_c^{-1} h_{c_1})_\# \in
C(PMod_R) = \{1\}$.\\

By using our sequence, with an inductive argument we get that
$(h_c)_\# = (h_v)_\#$ on $\mathcal{C}(R)$ and $(h_c)_\# = (h_v)_\#
= \lambda$ on $\{[c]\} \cup \{[v]\} \cup L_c \cup D_c \cup L_v
\cup D_v$. In particular we see that, $(h_c)_\#$ agrees with
$\lambda$ on any nonseparating curve $v$ and on $L_v$. Since every
separating curve is in the link, $L_r$, of some nonseparating
curve $r$, we see that $(h_c)_\#$ agrees with $\lambda$ on
$\mathcal{C}(R)$. This proves the lemma.\end{proof}\\

Proof of Theorem 1.1 follows from Lemma 3.13 and the remarks made
before this lemma. Note that in \cite{Ir}, we referred to Ivanov's
cutting arguments for the result of the corresponding theorem. The
proof of the
cutting argument given above can be adapted easily to that case.\\

\section{Injective Homomorphisms of Subgroups of Mapping Class Groups}

We assume that $\Gamma'=ker(\varphi)$ where $\varphi: Mod_R^*
\rightarrow Aut(H_1(R, \mathbb{Z}_3))$ is the homomorphism defined
by the action of homeomorphisms on the homology.

A mapping class $g \in Mod_R^*$ is called \textit{pseudo-Anosov}
if $\mathcal{A}$ is nonempty and if $g ^n (\alpha) \neq \alpha$,
for all $\alpha$ in $\mathcal{A}$ and any $n \neq 0$. $g$ is
called \textit{reducible} if there is a nonempty subset $
\mathcal{B} \subseteq \mathcal{A}$ such that a set of disjoint
representatives can be chosen for $\mathcal{B}$ and
$g(\mathcal{B}) = \mathcal{B}$. In this case, $ \mathcal{B}$ is
called a \textit{reduction system} for $g$. Each element of
$\mathcal{B}$ is called a \textit{reduction class} for $g$. A
reduction class, $\alpha$, for $g$, is called an \textit{essential
reduction class} for $g$, if for each $\beta \in \mathcal{A}$ such
that $i(\alpha, \beta) \neq 0$ and for each integer $m \neq 0$,
$g^m (\beta) \neq \beta$. The set, $\mathcal{B}_g$, of all
essential reduction classes for $g$ is called the
\textit{canonical reduction system} for $g$. The correspondence $g
\rightarrow \mathcal{B}_g$ is canonical. In particular, it
satisfies $g(\mathcal{B}_h) = \mathcal{B}_{ghg^{-1}}$ for all $g,
h$ in $Mod_R^*$.\\

The following two lemmas are well known facts. The isotopy class
of a Dehn twist about a circle $a$, is denoted by $t_{\alpha}$,
where $[a] = \alpha$.

\begin{lemma}
\label{wellknown1} Let $\alpha, \beta \in \mathcal{A}$ and $i, j$
be nonzero integers. Then, $t_\alpha ^i = t_\beta ^j
\Leftrightarrow \alpha = \beta$ and $i=j$.
\end{lemma}

\begin{lemma}
\label{wellknown2} Let $\alpha, \beta$ be distinct elements in
$\mathcal{A}$. Let $i, j$ be two nonzero integers. Then, $t_\alpha
^i t_\beta ^j = t_\beta ^j t_\alpha ^i \Leftrightarrow i(\alpha,
\beta)= 0$.
\end{lemma}

The proofs of the following two lemmas follow by the techniques
given in \cite{Ir}. Note that we need to use that the maximal rank
of an abelian subgroup of $Mod_R^*$ is $3g-3+p$, \cite{BLM}, in
these proofs.

\begin{lemma}
\label{rank=1} Let $K$ be a finite index subgroup of $Mod_R^*$ and
$f:K \rightarrow Mod_R^*$ be an injective homomorphism. Let
$\alpha \in \mathcal{A}$. Then there exists $N \in \mathbb{Z^*}$
such that

\begin{center}
 $rank$ $C(C_{\Gamma'} (f(t_{\alpha} ^{N})) ) = 1$.
\end{center}
\end{lemma}

\begin{lemma}
\label{reducible} Let K be a finite index subgroup of $Mod_R^*$.
Let $f:K \rightarrow Mod_R^*$ be an injective homomorphism. Then
there exists $N \in \mathbb{Z^*}$ such that $f(t_{\alpha} ^ N)$ is
a reducible element of infinite order for all $\alpha \in
\mathcal{A}$.
\end{lemma}

In the proof of Lemma \ref{reducible}, we use that centralizer of
a p-Anosov element in the extended mapping class group is a
virtually infinite cyclic group, \cite{Mc}.

\begin{lemma}
\label{correspondence} Let $K$ be a finite index subgroup of
$Mod_R^*$ and $f:K \rightarrow Mod_R^*$ be an injective
homomorphism. Then $\forall \alpha \in \mathcal{A}$, $f( t_\alpha
^N)= t_{\beta(\alpha)}^M$ for some $M, N \in \mathbb{Z^*}$,
$\beta(\alpha) \in \mathcal{A}$.
\end{lemma}

\begin{proof} Let $\Gamma= f^{-1}(\Gamma') \cap \Gamma'$. Since
$\Gamma$ is a finite index subgroup we can choose $N \in Z^*$ such
that $t_\alpha^N \in \Gamma$ for all $\alpha$ in $\mathcal{A}$. By
Lemma \ref{reducible} $f(t_{\alpha} ^ N)$ is a reducible element
of infinite order in $Mod_R^*$. Let $C$ be a realization of the
canonical reduction system of $f(t_{\alpha}^N)$. Let $c$ be the
number of components of $C$ and $r$ be the number of p-Anosov
components of $f(t_{\alpha} ^N)$. Since $t_{\alpha} ^ N \in
\Gamma, f(t_{\alpha} ^ N) \in \Gamma'$. By Theorem 5.9 \cite{IMc},
$C(C_{\Gamma'} (f(t_{\alpha} ^ N )))$ is a free abelian group of
rank $c+r$. By Lemma \ref{rank=1} $c+r=1$. Then, either $c=1$,
$r=0$ or $c=0$, $r=1$. Since there is at least one curve in the
canonical reduction system we have $c=1$, $r=0$. Hence, since
$f(t_{\alpha} ^ N) \in \Gamma'$, $f(t_{\alpha} ^{N}) = t_{\beta
({\alpha})}^{M}$ for some $M \in \mathbb{Z^*}$, $\beta(\alpha) \in
\mathcal{A}$, [1], [5].\end{proof}\\

\noindent {\bf Remark:} Suppose that $f(t_{\alpha} ^{M}) = t_\beta
^P$ for some $\beta \in \mathcal{A}$ and $M, P \in \mathbb{Z^*}$
and $f(t_{\alpha} ^{N}) = t_\gamma ^Q$ for some $\gamma \in
\mathcal{A}$ and $N, Q \in \mathbb{Z^*}$. Since $f(t_{\alpha} ^{M
\cdot N}) = f(t_{\alpha} ^{N \cdot M})$, $t_\beta ^{PN} = t_\gamma
^{QM}$, $P, Q, M, N \in \mathbb{Z^*}$. Then, $\beta = \gamma$ by
Lemma \ref{wellknown1}. Therefore, by Lemma \ref{correspondence},
$f$ gives a correspondence between isotopy classes of circles and
$f$ induces a map, $f_*: \mathcal{A} \rightarrow \mathcal{A}$,
where $f_*(\alpha) = \beta(\alpha)$.\\

In the following lemma we use a well known fact that $f t_\alpha
f^{-1}=t_{f(\alpha)} ^{\epsilon(f)}$ for all $\alpha$ in
$\mathcal{A}$, $f \in Mod_R^*$, where $\epsilon(f) = 1$ if $f$ has
an orientation preserving representative and $\epsilon(f) = -1$ if
$f$ has an orientation reversing representative.

\begin{lemma}
\label{identity} Let $K$ be a finite index subgroup of $Mod_R^*$.
Let $f:K \rightarrow Mod_R^*$ be an injective homomorphism. Assume
that there exists $N \in \mathbb{Z}^*$ such that $\forall \alpha
\in$ $\mathcal{A}$, $\exists Q \in \mathbb{Z}^*$ such that
$f(t_{\alpha} ^N) = t_{\alpha}^Q$. Then, $f$ is the identity on
$K$.
\end{lemma}

\begin{proof}
We use Ivanov's trick to see that $f(kt_{\alpha} ^ N k^{-1})=$
$f(t_{k(\alpha)} ^{\epsilon{(k)} \cdot N }) = t_{k(\alpha)} ^{Q
\cdot \epsilon{(k)}}$ and $f(kt_{\alpha} ^ N k^{-1}) = f(k)
f(t_{\alpha} ^N) f(k)^{-1}=$ $f(k) t_{\alpha} ^Q f(k)^{-1} =
t_{f(k)(\alpha)} ^{\epsilon(f(k))\cdot Q}$ $\forall \alpha \in
\mathcal{A}$, $\forall k \in K$. Then, we have $t_{k(\alpha)} ^{Q
\cdot \epsilon{(k)}} = t_{f(k)(\alpha)} ^{\epsilon(f(k)) \cdot Q}$
$\forall \alpha \in \mathcal{A}$, $\forall k \in K$. Hence,
$k(\alpha) = f(k)(\alpha)$ $\forall \alpha \in \mathcal{A}$,
$\forall k \in K$ by Lemma \ref{wellknown1}. Then,
$k^{-1}f(k)(\alpha) = \alpha$ $\forall \alpha \in \mathcal{A}$,
$\forall k \in K$. Therefore, $k^{-1}f(k)$ commutes with
$t_{\alpha}$ $\forall \alpha \in \mathcal{A}$, $\forall k \in K$.
Since $PMod_R$ is generated by Dehn twists, $k^{-1}f(k) \in
C(PMod_R)$ $\forall k \in K$. Since the genus of $R$ is at least
two and $R$ is not a closed surface of genus two, $C(PMod_R)$ is
trivial by 5.3 in [6]. So, $k = f(k)$ $\forall k \in K$. Hence,
$f=id_K$.\end{proof}

\begin{coroll}
\label{id} Let $g: Mod_R^* \rightarrow Mod_R^*$ be an isomorphism
and $h : Mod_R^* \rightarrow Mod_R^*$ be an injective
homomorphism.  Assume that there exists $N \in \mathbb{Z}^*$ such
that $\forall \alpha \in$ $\mathcal{A}$, $\exists Q \in
\mathbb{Z}^*$ such that $h(t_{\alpha} ^N) = g(t_{\alpha}^Q)$. Then
$g=h$.
\end{coroll}

\begin{proof} Apply Lemma \ref{identity} to $g^{-1} h$ with $K = Mod_R^*$.
Since for all $\alpha$ in $\mathcal{A}$, $g^{-1} h(t_{\alpha} ^N)
= t_{\alpha} ^{Q}$, we have $g^{-1} h = id_K$. Hence, $g =
h$.\end{proof}\\

By the remark after Lemma \ref{correspondence}, we have that $f: K
\rightarrow Mod_R^*$ induces a map $f_*: \mathcal{A} \rightarrow
\mathcal{A}$, where $K$ is a finite index subgroup of $Mod_R^*$.
In the following lemma we prove that $f_*$ is a superinjective
simplicial map on $\mathcal{C}(R)$.

\begin{lemma}
\label{intersection0} Let $f:K \rightarrow Mod_R^*$ be an
injection. Let $\alpha$, $\beta \in \mathcal{A}$. Then,
\[
i(\alpha,\beta)=0 \Leftrightarrow i(f_{*}(\alpha),
f_{*}(\beta))=0.
\]
\end{lemma}

\begin{proof} There exists $N \in \mathbb {Z^*}$ such that
$t_{\alpha} ^N \in K$ and $t_{\beta} ^N \in K$. Then we have the
following: $i(\alpha, \beta)=0$ $\Leftrightarrow$ $t_{\alpha} ^N
t_\beta ^N = t_{\beta} ^N t_\alpha ^N$ (by Lemma \ref{wellknown2})
$\Leftrightarrow$ $f(t_{\alpha} ^N) f(t_{\beta} ^ N) = f(t_{\beta}
^ N) f(t_{\alpha} ^ N)$ (since $f$ is injective on K)
$\Leftrightarrow$ $t_{f_*(\alpha)} ^P t_{f_*(\beta)}^Q =
t_{f_*(\beta)}^Q t_{f_*(\alpha)} ^P$ where $P = M(\alpha, N), Q =
M(\beta, N) \in \mathbb{Z}^*$ (by Lemma \ref{correspondence}) $
\Leftrightarrow i(f_{*}(\alpha),
f_{*}(\beta))=0$ (by Lemma \ref{wellknown2}).\end{proof}\\

Now, we prove the second main theorem of the paper.

\begin{theorem}
\label{main} Let $f$ be an injective homomorphism, $f:K
\rightarrow Mod_R^*$, then $f$ is induced by a homeomorphism of
the surface $R$ and $f$ has a unique extension to an automorphism
of $Mod_R^*$. \end{theorem}

\begin{proof} By Lemma \ref{intersection0} $f_*$ is a superinjective
simplicial map on $\mathcal{C}(R)$. Then, by Theorem 1.1, $f_*$ is
induced by a homeomorphism $h:R \rightarrow R$, i.e. $f_*(\alpha)
= h_\#(\alpha)$ for all $\alpha$ in $\mathcal{A}$, where
$h_\#=[h]$. Let $\chi ^ {h\#}: Mod_R^* \rightarrow Mod_R^*$ be the
isomorphism defined by the rule $\chi ^ {h_\#}(k) =
h_\#k{h_\#}^{-1}$ for all $k$ in $Mod_R^*$. Then for all $\alpha$
in $\mathcal{A}$, we have the following:

$\chi ^{h_\# ^{-1}} \circ f ({t_ \alpha} ^N) =  \chi ^{h_\# ^{-1}}
(t_{f_*(\alpha)}^ M) = \chi ^{h_\# ^{-1}} (t_{h_\#(\alpha)} ^M) =
h_\#^{-1} t_{h_\#(\alpha)} ^M h_\# = t_ {h_\#^{-1} (h_\#(\alpha))}
^{M \cdot \epsilon{(h_\#^{-1})}} = t_\alpha ^{M \cdot
\epsilon{(h_\#^{-1})}}$.

Then, since $\chi ^{h_\#^{-1}} \circ f$ is injective, $\chi
^{h_\#^{-1}} \circ f = id_K$ by Lemma \ref{identity}. So, $\chi
^h_\# |_K = f$. Hence, $f$ is the restriction of an isomorphism
which is conjugation by $h_\#$, (i.e. $f$ is induced by $h$).

Suppose that there exists an automorphism $\tau :  Mod_R^*
\rightarrow Mod_R^*$ such that $\tau |_{K}=f$. Let $N \in Z^*$
such that $ t_\alpha ^N \in K$  for all $\alpha$ in $\mathcal{A}$.
Since $\chi ^h_\# |_K = f = \tau |_K$ and $t_\alpha ^N \in K$,
$\tau(t_\alpha ^N) = \chi ^h_\#(t_\alpha ^N)$ for all $\alpha$ in
$\mathcal{A}$. Then, by Corollary \ref{id}, $\tau = \chi ^
{h_\#}$. Hence, the extension of $f$ is unique. \end{proof}\\

{\bf Acknowledgments}\\

We thank John D. McCarthy for his suggestions and comments on this
paper. We also thank Nikolai Ivanov for his comments and Ursula
Hamenstaedt for helpful discussions.

\vspace{1cm}

\noindent University of Michigan, Department of Mathematics, Ann
Arbor, MI 48109, USA; eirmak@umich.edu\\

\end{document}